\documentclass[10pt,a4paper,notitlepage,reqno]{amsart}

\usepackage[utf8]{inputenc}                          
\usepackage[english]{babel}                          

\usepackage[margin=1in]{geometry}                    
\usepackage{amsmath,amssymb,amsthm,amsfonts,mathrsfs} 
\usepackage{bbm}                                      
\usepackage{dsfont}                                   
\usepackage{latexsym}                                 
\usepackage{wasysym}                                  
\usepackage{arydshln}                                 

\usepackage[usenames,dvipsnames]{xcolor}              
\usepackage{setspace}                                 
\usepackage{enumitem}                                 
\usepackage[font=footnotesize]{caption}               
\usepackage{subcaption}                               
\usepackage{booktabs,hhline,multirow}                 
\usepackage{multicol}                                 

\usepackage{hyperref}                                 
\usepackage{url}                                      

\usepackage{graphicx}                                 

\usepackage{tikz}                                     
\usetikzlibrary{matrix,chains,shapes,arrows,arrows.meta,scopes,positioning,
  calc,decorations.markings,decorations.pathreplacing,decorations.pathmorphing,patterns}

\tikzstyle{empty}  =[circle,draw=black!80,thick]      
\tikzstyle{emptyn} =[circle,draw=black!80,fill=white,scale=0.5]
\tikzstyle{nero}   =[circle,draw=black!80,fill=black!80,thick]

\usepackage[vcentermath]{youngtab}                    

\allowdisplaybreaks                                   
\newtheorem{theorem}{Theorem}[section]                
\newtheorem*{theorem*}{Theorem}                       
\newtheorem{lemma}[theorem]{Lemma}                    
\newtheorem{proposition}[theorem]{Proposition}        
\newtheorem*{proposition*}{Proposition}               
\newtheorem{corollary}[theorem]{Corollary}            
\newtheorem*{corollary*}{Corollary}                   

\theoremstyle{definition}                             
\newtheorem{definition}[theorem]{Definition}          
\newtheorem{notation}[theorem]{Notation}              
\newtheorem{remark}[theorem]{Remark}                  
\newtheorem{example}[theorem]{Example}                

\numberwithin{equation}{section}                      

\newcommand{\Fix}{\operatorname{Fix}}
\newcommand{\Char}{\operatorname{Char}}               
\newcommand{\Gal}{\operatorname{Gal}}                 
\newcommand{\Irr}{\operatorname{Irr}}                 
\newcommand{\Syl}{\operatorname{Syl}}                 
\newcommand{\Lin}{\operatorname{Lin}}                 

\newcommand{\C}{\mathbb{C}}                           
\newcommand{\Q}{\mathbb{Q}}                           
\newcommand{\Z}{\mathbb{Z}}                           
\newcommand{\N}{\mathbb{N}}                           
\newcommand{\triv}{\mathbbm{1}}                       
\newcommand{\Aut}{\mathrm{Aut}}                        
                        
\newcommand{\up}{\big\uparrow} 
\newcommand{\down}{\big\downarrow} 
\newcommand{\giu}{\downarrow} 

\newcommand{\cF}{\mathcal{F}} 
\newcommand{\cG}{\mathcal{G}}

\newcommand{\cT}{\mathcal{T}} 
\newcommand{\cX}{\mathcal{X}}

\newcommand{\cR}{\mathcal{R}} 

\newcommand{\sk}{\mathfrak{s}} 
\newcommand{\ts}{\mathtt{s}} 
\newcommand{\tF}{\mathtt{F}} 
\newcommand{\bF}{\mathbf{F}} 

\begin{document}

\title{On the local representation theory of symmetric groups}  


\author{Greta Tendi}
\address[G. Tendi]{Università di Pisa, Dipartimento di Matematica, Largo Bruno Pontecorvo 5, Pisa, Italy}
\email{greta.tendi@phd.unipi.it}


\begin{abstract}
Let $p$ be a prime. Given a Sylow $p$-subgroup $P$ of a finite symmetric group, we describe the action of its normalizer on $\mathrm{Irr}(P)$. To this end, we establish a one-to-one correspondence between the irreducible characters of $P$ and certain equivalence classes of explicitly defined functions, which are also naturally suited to describing the Galois action.
\end{abstract}

\maketitle




\section{Introduction}\label{sec:intro}

Partially motivated by the study of the McKay conjecture and its refinements, recent years have seen increased attention to the relationship between the irreducible characters of the symmetric group $S_n$ and those of its Sylow $p$-subgroups. For instance, a complete description of the action of the normalizer subgroup $N_{S_n}(P)$ on the set of linear characters $\mathrm{Lin}(P)$ is given in \cite{GANT}, where $P$ denotes a Sylow \mbox{$p$-subgroup} of $S_n$. This result is then used as a key ingredient to show that the symmetric and alternating groups satisfy a strengthening of the McKay-Navarro conjecture (see \cite[Theorem A]{GANT}), providing an alternative proof to the same theorem previously established by Navarro and Tiep \cite{NT}.
On the other hand, in \cite{GL3} the authors describe the set of irreducible constituents arising from the induction of any irreducible character of $P$ to $S_n$. A fundamental tool in obtaining this result is the parametrization of the set of irreducible characters $\mathrm{Irr}(P)$ by certain combinatorial objects, namely tuples of rooted, complete, $p$-ary, labeled trees.

In this article, we replace the previous parametrization with a new function-based one, which is conceptually equivalent but offers a clearer and more explicit framework for the analysis of $\mathrm{Irr}(P)$. In particular, we employ this reformulation to completely determine the action of the normalizer subgroup $N_{S_n}(P)$ on the set $\mathrm{Irr}(P)$. This constitutes the first main result of this work, it is established in \textbf{Theorems~\ref{th:1} and~\ref{th:2}} and significantly generalizes the preceding description of the action on linear characters given in \cite{GANT}.
The second main result, contained in \textbf{Theorems~\ref{th:3} and~\ref{th:4}}, provides a precise combinatorial description of the action of the Galois group on $\mathrm{Irr}(P)$. Furthermore, \mbox{\textbf{Theorem~\ref{th:5}}} yields a straightforward method to determine the product of a character in $\Irr(P)$ with one in $\Lin(P)$, illustrating the efficiency of this formalism for computational purposes.

We remark that in \cite{law} it was observed that if two linear characters of $P$ are Galois-conjugate, then they also lie in the same orbit under the action of $N_{S_n}(P)$. Moreover, in \cite{gll} it was proved that two linear characters of $P$ are $N_{S_n}(P)$-conjugate if and only if their induction from $P$ to $S_n$ yields the same character. As a consequence of Theorems~\ref{th:2} and~\ref{th:4}, we show that these statements do not extend to the full set of irreducible characters $\mathrm{Irr}(P)$; explicit counterexamples are provided in \textbf{Proposition~\ref{rem:conclusions}}.

The paper is structured as follows. We begin by recalling some basic representation-theoretic properties of Sylow subgroups of symmetric groups. After introducing the necessary notation, in Section~\ref{sec:chT-functions} we index $\mathrm{Irr}(P)$ by \textit{$\cT$-functions} (see Definition~\ref{def:admT}), thereby establishing a bijection in a setting shown to be equivalent to the tree-based model presented in \cite{GL3}. In Section~\ref{sec:norgalIrr(P)}, we use this parametrization to describe the actions of the normalizer and of the Galois group on these characters, as well as to compute the product of an irreducible character with a linear one. Finally, in Section~\ref{sec:eqonch}, we discuss further consequences of our results, with a focus on their connection to character induction (see Table~\ref{tab:schema}). 

\subsection*{Acknowledgments}
The author thanks Eugenio Giannelli for suggesting and supervising this work, and the anonymous referee for helpful comments and valuable suggestions on a previous version of this paper.

\section{Sylow subgroups and Sylow normalizers of symmetric groups}\label{sec:prelim}

In this section we fix the notation and summarize the main representation-theoretic facts required later. 
We refer the reader to \cite{Is, Nav98, JK} for a comprehensive account of these topics.

\begin{notation}
For any pair of integers $x,y\in\Z$ we write $[x,y]$ for the set $\{z\in\Z\ |\ x\leq z\leq y\}$. Given any $n\in\N$, we denote by $S_n$ the symmetric group acting on a set of cardinality $n$. Moreover, we will say that a finite collection of objects is \emph{mixed} if the objects are not all equal (that is, if the corresponding set has more than one element).
\end{notation}

\subsection{Representations of wreath products}\label{sec:repofwp}
Consider $G$ as a finite group and $H\le G$ a subgroup. We write $\Char(G)$ for the set of ordinary characters of $G$, and $\Irr(G)$ and $\Lin(G)$ for the subsets of irreducible and linear characters, respectively. For $\chi\in\Irr(G)$ and $\varphi\in\Irr(H)$, we denote by $\chi\down_H$ the restriction of $\chi$ to $H$, and by $\varphi\up^G$ the induction of $\varphi$ to $G$. We use square brackets $[\cdot,\cdot]$ to denote the standard inner product of characters and let $\Irr(G\mid\varphi):=\{\chi\in\Irr(G) \mid [\chi\down_H,\varphi] \ne 0\}$ be the subset of all the irreducible constituents of the induced character $\varphi\up^G$.

We now fix our notation for characters of wreath products, adopting the same conventions as in \cite[Chapter 4]{JK}. Given a natural number $n$, we denote by $G^{\times n}$ the direct product of $n$ copies of $G$. For any subgroup $H\le S_n$, the permutation action of $S_n$ on the direct factors of $G^{\times n}$ induces an action of $S_n$ (and therefore of $H\le S_n$) via automorphisms on $G^{\times n}$, giving the wreath product $G\wr H:= G^{\times n}\rtimes H$. The normal subgroup $G^{\times n}$ is sometimes called the base group of the wreath product (whereas $H$ is the top group), and the elements of $G\wr H$ are denoted by $(g;h):=(g_1,\dotsc,g_n;h)$ for $g_i\in G$ and $h\in H$. 

Let $V$ be a $\C G$--module and let $\phi$ be the character afforded by $V$. We let $V^{\otimes n}:=V\otimes\cdots\otimes V$ ($n$ copies) be the corresponding $\C G^{\times n}$--module. The left action of $G\wr H$ on $V^{\otimes n}$ defined by linearly extending
$$(g_1,\dotsc,g_n;h)\ :\quad v_1\otimes \cdots\otimes v_n \longmapsto g_1v_{h^{-1}(1)}\otimes\cdots\otimes g_nv_{h^{-1}(n)}$$
turns $V^{\otimes n}$ into a $\C (G\wr H)$--module, which we denote by $\widetilde{V^{\otimes n}}$. We let $\widetilde{\phi^{\times n}}$ denote the character afforded by the representation $\widetilde{V^{\otimes n}}$. For any character $\psi$ of $H$, we abuse notation and let $\psi$ also denote its inflation to $G\wr H$. Finally, we introduce the symbol
\[
\cX(\phi; \psi) := \widetilde{\phi^{\times n}} \cdot \psi \in\Char(G\wr H)
\]
to denote the  character of $G\wr H$ obtained as the product of $\widetilde{\phi^{\times n}}$ and $\psi$.

Recall that the irreducible characters of $G^{\times n}$ are precisely the direct products $\zeta_1\times\cdots\times\zeta_n$, where $\zeta_i\in\Irr(G)$. Let $\phi\in\Irr(G)$ and set $\phi^{\times n}:= \phi \times \cdots \times \phi\in\Irr(G^{\times n})$. Observe that $\widetilde{\phi^{\times n}}\in\Irr(G\wr H)$ is an extension of $\phi^{\times n}$. Hence, by Gallagher's Theorem \cite[Corollary 6.17]{Is} we have
\[
\Irr(G\wr H \mid \phi^{\times n})= \{ \cX(\phi; \psi) \mid \psi \in \Irr(H)\}.
\]

Now, let $p$ be a prime and let $H=C_p$ denote the cyclic group of order $p$. It is well known that Sylow $p$-subgroups of symmetric groups are isomorphic to direct products of copies of iterated wreath products of $C_p$. To better understand the representation theory of these groups, we now revisit the preceding discussion, focusing on this specific case. By basic Clifford theory, we obtain:

\begin{lemma}\label{lem:clif}
Each $\eta \in \Irr(G\wr C_p)$ falls into exactly one of the following two cases.
\begin{itemize}
	\item[(i)] $\eta=\zeta_1\times\cdots\times\zeta_p\up^{G\wr C_p}$ for some mixed $\zeta_1,\dotsc,\zeta_p\in\Irr(G)$, so that $\Irr(G\wr C_p\mid \zeta_1\times\cdots\times\zeta_p)=\{\eta\}$ and $\eta\down_{G^{\times p}}=\sum_{\substack{\pi\in C_p}}\zeta_{(1)\pi}\times\cdots\times\zeta_{(p)\pi}$.	
	\item[(ii)] $\eta=\cX(\phi;\psi)$ for some $\phi\in\Irr(G)$ and $\psi\in\Irr(C_p)$, so that $\phi^{\times p}\up^{G\wr C_p}=\sum_{\substack{\varphi\in\Irr(C_p)}} \cX(\phi;\varphi)$	and	$\eta\down_{G^{\times p}}=\phi^{\times p}$.
\end{itemize}
\end{lemma}

\subsection{Characters of Sylow subgroups of symmetric groups}\label{sec:chofSylow}
Given any prime $p$ and $n\in\N$, let $P_n\in\Syl_p(S_n)$ be a Sylow $p$-subgroup of $S_n$.
We observe that $P_1$ is the trivial group, and that $P_p\cong C_p$. For every integer $k\ge 2$, we have that 
\[
P_{p^k}\cong\big(P_{p^{k-1}}\big)^{\times p} \rtimes P_p=P_{p^{k-1}}\wr P_p\cong P_p\wr \cdots \wr P_p
\]
that is, a $k$-fold wreath product. In general, if the $p$-adic expansion of $n$ is 
\[
n=\sum_{i=1}^{t}a_ip^{k_i}
\]
where $t\in\N$, $k_1>\cdots>k_t\geq 0$ are integers, and $a_1,\dotsc,a_t\in [1,p-1]$, then 
\begin{equation}\label{eq:sylowprod}
P_n\cong (P_{p^{k_1}})^{\times a_1}\times\cdots\times (P_{p^{k_t}})^{\times a_t}.
\end{equation}
We derive that, setting  
\begin{equation}\label{eq:pgencar}
\begin{split}
&q_0 :=0 \quad \text{and } \quad q_{i}:=q_{i-1}+a_{i} \quad \text{for all } i\in [1,t], \\
\Irr(P_n) &= \{\theta_1\times\cdots\times\theta_{q_t}\ |\  \theta_j\in\Irr(P_{p^{k_i}}) \ \text{for}\ i\in [1,t], \ j\in[q_{i-1}+1,q_i] \}.
\end{split}
\end{equation}
Thus, it suffices to study $\Irr(P_{p^k})$ for $k\ge 0$ to describe the irreducible characters of any $P_n$. We note that $P_1$ admits only the trivial character, while $\Irr(P_p)=\Lin(P_p)$ consists of $p$ characters of degree one, as $P_p$ is abelian. If $k>1$, every irreducible character of $P_{p^k}$ can be expressed recursively in terms of those of $P_p$, as recorded in the next lemma. This result follows immediately from iterated applications of Lemma~\ref{lem:clif}, which describes the irreducible characters of wreath products with top group $C_p$. 

\begin{lemma}\label{lem:chP-trees}
Let $p$ be a prime. For every integer $k\geq 1$ and $P_{p^k}\in\Syl_p(S_{p^k})$, if $\theta\in\Irr(P_{p^k})$ then exactly one of the following two cases holds:
\begin{itemize}
	\item[(i)] $\theta= \theta_1 \times \cdots \times \theta_p \up^{P_{p^k}}$, for some mixed $\theta_1 , \dotsc, \theta_p \in \Irr(P_{p^{k-1}})$, or
	\item[(ii)] $\theta=\cX(\phi;\psi)$ for some $\phi\in\Irr(P_{p^{k-1}})$ and $\psi\in \Irr(P_p)$. 
\end{itemize}
\end{lemma}

\begin{remark}\label{rem:chlimits}
Notice that in case (i) of the lemma above, $\theta=\theta_{(1)\pi}\times\cdots\times\theta_{(p)\pi}\up^{P_{p^k}}$ for any cyclic permutation $\pi\in P_p$. 
\end{remark}

\subsection{Action of the normalizer of a Sylow subgroup in a symmetric group}\label{sec:gensylnor}
We conclude this preliminary section by recalling explicit presentations of a fixed Sylow $p$-subgroup $P_n$ of $S_n$ and of its normalizer $N_n$. Following \cite[Section 2]{GANT}, we recover generators for both groups and describe their interaction under conjugation. We then extend the results therein by characterizing the inherited action of $N_n/P_n$ on $P_n$. This will be used in Section~\ref{sec:norgalIrr(P)} to study the induced action of $N_n$ on the set $\Irr(P_n)$. Examples~\ref{ex:111} and~\ref{ex222} serve as a guide throughout. 

\subsubsection{Prime power degree case}\label{sec:gensylnorppc}
Every Sylow $p$-subgroup of $S_{p^k}$ is isomorphic to the $k$-fold wreath product $C_p\wr \dots \wr C_p$, and hence admits a minimal generating set of cardinality $k$. We construct such a set by fixing the nontrivial $p$-cycle $\gamma:=(1,2,\dotsc,p) \in C_p \le S_p$, and, at each step in the recursive construction of the wreath product, selecting the element in the top group that permutes the $p$ factors of the base group according to $\gamma$. Define $\gamma_1:=\gamma\in P_p$, and for $j\in [2,k]$, let $\gamma_j:=(1,1,\dotsc, 1;\gamma)\in P_{p^{j-1}}\wr P_p=P_{p^j}$.
Let $\psi_j : P_{p^j} \to S_{p^j}$ be the canonical permutation representation of the $j$-fold wreath product (see \cite[4.1.18, Chapter 4]{JK}), and let $\iota_j : S_{p^j} \hookrightarrow S_{p^k}$ be the natural embedding. Set $f_j:=\iota_j\circ\psi_j$ and define $g_j^{(k)}:=f_j(\gamma_j)$ for each $j\in [1,k]$. The resulting $k$ permutations are elements of $S_{p^k}$ and act on the set $[1,p^k]$. For clarity, we include in the following definition their disjoint cycle decomposition. 

\begin{definition}\label{def:P_p^k}
Let $k\in \N$. For each $i\in [1,k]$,
\[
g_i^{(k)} = \prod_{j=1}^{p^{i-1}} (j\ ,\ p^{i-1}+j\ ,\ 2p^{i-1}+j\ ,\ \cdots\ ,\ (p-1)p^{i-1}+j) \in S_{p^k}.
\]
We denote by $P_{p^k}$ the specific Sylow $p$-subgroup of $S_{p^k}$ generated by $\{g_i^{(k)}\mid i\in[1,k]\}$:
\[
P_{p^k} := \langle g_1^{(k)}, \dotsc, g_k^{(k)} \rangle.
\]
\end{definition}

\begin{remark}\label{rem:obsgen}
Observe that for each $j<k$, the restriction of $g_j^{(k)}$ to $[1,p^j]$ is $g_j^{(j)}$, and these permutations have the same disjoint cycle decomposition up to fixed points. With a slight abuse of notation, we have
\[
g_k^{(k)}=\gamma_k,\ \ g_{k-1}^{(k)} = (\gamma_{k-1},1,\dotsc,1;1),\ \ g_{k-2}^{(k)}=\big( (\gamma_{k-2},1,\dotsc,1;1), 1,\dotsc,1; 1\big)
\] 
and so on. In general, for all $j\in [1,k-1]$,  
\[
g_j^{(k)}=\big(\big(\dotsc \big( (\gamma_j,1,\dotsc,1;1),1,\dotsc,1; 1\big)\dotsc\big),1,\dotsc,1; 1\big). 
\]
Equivalently, by associativity of wreath products and compatibility with the canonical permutation representations, for all $j<k$ the generators satisfy the recursive relation 
\[
g_j^{(k)}=(g_j^{(k-1)}, 1, \dotsc, 1; 1)\in P_{p^{k-1}}\wr P_p=P_{p^{k}}.
\]
\end{remark}

The normalizer has the well-known structure $$N_{S_{p^k}}(P_{p^k})\cong P_{p^k}\rtimes (C_{p-1})^{\times k}$$ (see, for instance, \cite[Lemma 4.1]{OlssonPaper}), and admits a minimal generating set $\{g_1^{(k)},\dotsc,g_k^{(k)},\sigma_1^{(k)},\dotsc,\sigma_k^{(k)} \}$ of cardinality $2k$. As with the generators of $P_{p^k}$, the remaining $k$ generators may also be described recursively. To introduce these elements, we require the following definition.

\begin{definition}
Let $p$ be a prime. We say that $a\in\Z$ is a primitive root modulo $p$ if its multiplicative order modulo $p$ is $p-1$, so that $\langle a+p\Z \rangle \cong (\Z/p\Z)^{\times}$.
\end{definition}

\begin{example}\label{ex:primroot}
Recall that $(\Z/p\Z)^{\times}$ is the multiplicative group of $\Z/p\Z$, and that we can always find a primitive root modulo $p$ in the set $[1,p-1]$. For example, if $p=7$ then $3$ is a primitive root modulo $7$, since its multiplicative order modulo $7$ is $6$: $3^1\equiv 3, \quad 3^2\equiv 2, \quad 3^3\equiv 6, \quad 3^4\equiv 4, \quad 3^5\equiv 5$ and $3^6\equiv 1\pmod{7}$. Note that $(1,2,3,4,5,6,7)^{(3,2,6,4,5,1)}=(3,6,2,5,1,4,7)=
(1,2,3,4,5,6,7)^3$.
\end{example}

We now define the permutations $\sigma_j^{(k)}\in S_{p^k}$ recursively, and then detail their characteristic properties.

\begin{definition}\label{def:sigmajk}
Let $a,a_1,\dotsc,a_{p-1}\in [1,p-1]$ where $a$ is a primitive root modulo $p$, and $a_i\equiv a^i\pmod{p}$ for all $i\in [1,p]$. For any $k\in \N$ and $m\in\Z$, let $\tau_m\in S_{p^k}$ be the permutation $i\mapsto i+m$ with numbers modulo $p^k$ (taken in the range $[1,p^k]$). For $1\le j<k$, set
\[
\sigma_j^{(k)} = \prod_{i=0}^{p-1} (\sigma_j^{(k-1)})^{\tau_{ip^{k-1}}} \text{ \ and \ } \sigma_k^{(k)} = \prod_{i=0}^{p^{k-1}-1} (a_1p^{k-1}, a_2p^{k-1}, \dotsc, a_{p-1}p^{k-1})^{\tau_{-i}}
\]
with numbers modulo $p^k$ (taken in the range $[1,p^k]$). For notational convenience, we set $\tau:=(\sigma_1^{(1)})^{-1}$.
\end{definition}

As shown in Definition~\ref{def:P_p^k}, for every $j\in [1,k]$ the permutation $g_j^{(k)}$ is a product of disjoint $p$-cycles. The permutations $\{(g_j^{(k)})^{\tau_{sp^j}}\mid s\in [0,p^{k-j}-1]\}$ are products of disjoint $p$-cycles. Let $C$ be the set of all the $p$-cycles appearing in the decompositions of the $(g_j^{(k)})^{\tau_{sp^j}}$, as $s$ varies; these cycles are mutually disjoint and partition $[1,p^k]$. The element $\sigma_j^{(k)}$ is simply the product, over all $p$-cycles $(c_1,c_2,\dotsc,c_p)$ in $C$, of $(p-1)$-cycles $(c_{a_1},c_{a_2},\dotsc,c_{a_{p-1}})$. It follows that, for each $p$-cycle $(c_1,c_2,\dotsc,c_p)$ in $C$, we have $\sigma_j^{(k)}(c_p)=c_p$, while $c_i\notin\Fix(\sigma_j^{(k)})$ for $i\neq p$. Equivalently, $\sigma_j^{(k)}$ is the unique permutation in $S_{p^k}$ such that 
\begin{equation}\label{eq:infosigmajk}
\begin{aligned}
\Fix(\sigma_j^{(k)}) &= \bigcup_{s=0}^{p^{k-j}-1} \{ ((p-1)p^{j-1}+i)\tau_{sp^j} \mid i \in [1,p^{j-1}] \} \\ 
((g_j^{(k)})^{\tau_{sp^j}})^{\sigma_j^{(k)}} &= ((g_j^{(k)})^{\tau_{sp^j}})^a \qquad \quad \text{for all } s\in [0,p^{k-j}-1].
\end{aligned}
\end{equation}
Thus, $\sigma_j^{(k)}$ is a product of disjoint $(p-1)$-cycles and normalizes $\langle (g_j^{(k)})^{\tau_{sp^j}} \rangle$ for every $s\in [0,p^{k-j}-1]$. In particular, the action of $\sigma_j^{(k)}$ on $[1,p^k]$ is uniquely determined by \eqref{eq:infosigmajk} for all $j\in [1,k]$. 

We now show that $\langle \sigma_t^{(k)} \mid t\in [1,k] \rangle \le N_{S_{p^k}}(P_{p^k})$, and that for each $j\in [1,k]$, the inner automorphism of $S_{p^k}$ induced by $\sigma_j^{(k)}$ restricts to an automorphism of $P_{p^k}$, whose action is explicitly described in Proposition~\ref{prop:sigmasuP}.

\begin{notation}
We record some useful notational remarks:
\begin{itemize}
	\item[(1)] We write $(P_{p^{k-1}})_m$ to denote the $m$th direct factor of the base group $B=(P_{p^{k-1}})^{\times p}$, where $P_{p^k}= B \rtimes P_p$.
	\item[(2)] For any permutation $\pi\in S_p$ and any $(h;\lambda)=(h_1,\dotsc,h_p;\lambda)\in P_{p^{k-1}} \wr P_p=P_{p^{k-1}}$, we write $(h^{\pi}, \lambda):=(h_{(1)(\pi)^{-1}}, \dotsc, h_{(p)(\pi)^{-1}}; \lambda)\in B\rtimes P_p=P_{p^k}$. 
	\item[(3)] We recall that $\tau=(\sigma_1^{(1)})^{-1}$ and that $(\gamma)^{\sigma_1^{(1)}}=\gamma^{a}$, where $\gamma=(1,2,\dotsc,p)$ and $a$ is the primitive root modulo $p$ fixed in Definition~\ref{def:sigmajk}.
\end{itemize}
\end{notation}

\begin{proposition}\label{prop:sigmasuP} 
Let $p$ be a prime and $1\le k\in\N$.
Let $(f, \delta)=(f_1, \dotsc, f_p;\delta)\in P_{p^{k-1}} \wr P_p=P_{p^k}$. Consider $\sigma_j^{(k)}$ for $j\in [1,k]$. Then,
\begin{equation*}
(f,\delta)^{\sigma_j^{(k)}}=
\begin{cases}
(f^{\sigma_1^{(1)}}; \delta^{\sigma_1^{(1)}})=(f_{(1)\tau}, \dotsc, f_{(p)\tau}; \delta^{\sigma_1^{(1)}}) & \text{if \ } j=k \\ 
(f_1^{\sigma_j^{(k-1)}},\dotsc,f_p^{\sigma_j^{(k-1)}}; \delta) & \text{if \ } j<k
\end{cases}
\end{equation*} 
\end{proposition}

\begin{proof}
We proceed by induction on $k$. If $k=1$, then  $j=k$, $\sigma_j^{(k)}=\sigma_1^{(1)}$ and $(\delta)^{\sigma_1^{(1)}}=(\delta)^a$. Only the first case arises and the proposition holds in this instance. Let $k>1$. We collect some immediate consequences of Definition~\ref{def:P_p^k}, Remark~\ref{rem:obsgen} and Definition~\ref{def:sigmajk}. 
Let $r,t,d,\ell \in [1,k-1]$, with $r<t$.
\begin{itemize} 
	\item[(i)] $(g_t^{(k)})^{g_k^{(k)}}=(g_t^{(k-1)}, 1, \dotsc, 1;1)^{g_k^{(k)}}=(1, g_t^{(k-1)}, \dotsc,1;1)$.
	\item[(ii)] $(g_r^{(k)})^{g_t^{(k)}} = (( \dots ((g_r^{(t)})^{g_t^{(t)}}) \dots ), 1, \dotsc, 1;1) = (g_r^{(k)})^{\tau_{p^{t-1}}}$.
	\item[(iii)] $(g_t^{(k)})^{\sigma_\ell^{(k)}}=(g_t^{(k-1)}, 1, \ldots, 1; 1)^{\sigma_\ell^{(k)}}=((g_t^{(k-1)})^{\sigma_\ell^{(k-1)}}, 1, \ldots, 1; 1)$ .
	\item[(iv)]$(g_k^{(k)})^{\sigma_d^{(k)}}= g_k^{(k)}$, and $\sigma_d^{(k)}$ permutes the cycles in the decomposition of $g_k^{(k)}$.
	\item[(v)] $(g_t^{(\ell)})^{\sigma_r^{(\ell)}}= g_t^{(\ell)}$ for all $\ell>t$. Indeed, proceeding by induction and using (iii)-(iv), if $\ell=t+1$, then $(g_t^{(t+1)})^{\sigma_r^{(t+1)}}=((g_t^{(t)})^{\sigma_r^{(t)}},1,\ldots, 1; 1)=(g_t^{(t)},1,\ldots, 1; 1)=g_t^{(t+1)}$. If $\ell>t+1$, then by the inductive hypothesis $(g_t^{(\ell)})^{\sigma_r^{(\ell)}}=((g_t^{(\ell-1)})^{\sigma_r^{(\ell-1)}},1,\ldots, 1; 1)=(g_t^{(\ell-1)},1,\ldots, 1; 1)=g_t^{(\ell)}$.
	\item[(vi)] $P_{p^{\ell}} \cong \langle g_t^{(k)} \mid t\in[1,\ell] \rangle$, which is a subgroup of $(P_{p^{k-1}})_1$. From now on we identify $P_{p^{\ell}}$ with its isomorphic copy lying in the first direct factor of the base group of $P_{p^k}$, for every $\ell<k$. This natural embedding induces an inclusion ordering among these subgroups, allowing us to omit superscripts and simply write $g_t$ instead of $g_t^{(k)}$ for $t\in[1,k]$.
\end{itemize}
For each $i\in[1,p]$, there exists an integer $n_i\in \N$ and a function $h_i : [1,n_i] \to [1,k-1]$ such that $f_i=g_{h_i(1)} \dots g_{h_i(n_i)}$. Let $h: \{ (i,j)\in \N^{\times 2} \mid i\in[1,p],  j\in[1,n_i] \} \to [1,k-1]$ be such that $h(i,j)=h_i(j)$ for any $(i,j)\in [1,p]\times[1,n_i]$, so that
\begin{equation}\label{eq:deffi}
f_i=g_{h(i,1)}\dots g_{h(i,n_i)}=\prod_{\ell=1}^{n_i} g_{h(i,\ell)}
\end{equation}
for all $i\in [1,p]$. If $\delta=\gamma^s$ for some $s\in [0,p-1]$, then
\begin{align*}
&(f;\delta)=(f_1,\dotsc,f_p;\delta)=(g_{h(1,1)} \dots g_{h(1,n_1)}, \dotsc, g_{h(p,1)} \dots g_{h(p,n_p)};\gamma^s) \\
&=(g_{h(1,1)} \dots g_{h(1,n_1)},1,\dotsc,1;1) \dots (1,\dotsc,1,g_{h(p,1)} \dots g_{h(p,n_p)};1)(1,\dotsc,1;\gamma)^s \\
&=(g_{h(1,1)},1,\dotsc,1;1)\dots(g_{h(1,n_1)},1,\dotsc,1;1)\dots(1,\dotsc,1,g_{h(p,1)};1)\dots(1,\dotsc,1,g_{h(p,n_p)};1)(1,\dotsc,1;\gamma)^s.
\end{align*}
Using (i), we can write $(f;\delta)$ solely in terms of the generators of $P_{p^k}$. Explicitly, we obtain
\begin{align}\label{eq:scritturagen}
(f;\delta)=g_{h(1,1)}\dots g_{h(1,n_1)} \dots (g_{h(p,1)})^{(g_k)^{p-1}} \dots (g_{h(p,n_p)})^{(g_k)^{p-1}}(g_k)^s=(\prod_{i=1}^{p}\prod_{\ell=1}^{n_i}(g_{h(i,\ell)})^{(g_k)^{i-1}})(g_k)^{s}.
\end{align}
We consider a representative factor of the product and express it using the disjoint cycle decomposition of the $g_t$. For $i\in[1,p]$ and $\ell\in[1,n_i]$, each conjugate element decomposes into disjoint $p$-cycles as
\begin{equation}\label{eq:scritturaprod}
\begin{aligned}
&(g_{h(i,\ell)})^{(g_k)^{i-1}}= \\
&=\prod_{t=1}^{p^{h(i,\ell)-1}}(t+(i-1)p^{k-1},t+p^{h(i,\ell)-1}+(i-1)p^{k-1},\dots,t+(p-1)p^{h(i,\ell)-1}+(i-1)p^{k-1}).
\end{aligned}
\end{equation}
In view of \eqref{eq:scritturagen} and \eqref{eq:scritturaprod}, we are now ready to compute
\begin{align}\label{eq:star}
(f;\delta)^{\sigma_j^{(k)}}&=(\prod_{i=1}^{p}\prod_{\ell=1}^{n_i}(g_{h(i,\ell)})^{(g_k)^{i-1}})^{\sigma_j^{(k)}}((g_k)^s)^{\sigma_j^{(k)}}=(\prod_{i=1}^{p}\prod_{\ell=1}^{n_i}(g_{h(i,\ell)})^{(g_k)^{i-1} \cdot \sigma_j^{(k)}})(g_k^{\sigma_j^{(k)}})^s.
\end{align}
We distinguish two cases. 

\noindent\textbf{Case j=k.} 
By \eqref{eq:scritturaprod} and the definition of $\sigma_k^{(k)}$, for any $i\in[1,p]$ and $\ell\in[1,n_i]$ one has:
\begin{align*}
&((g_{h(i,\ell)})^{(g_k)^{i-1}})^{\sigma_k^{(k)}}= \\
&=\prod_{t=1}^{p^{h(i,\ell)-1}}(t+(i-1)p^{k-1},t+p^{h(i,\ell)-1}+(i-1)p^{k-1},\dotsc,t+(p-1)p^{h(i,\ell)-1}+(i-1)p^{k-1})^{\sigma_k^{(k)}} \\
&=\prod_{t=1}^{p^{h(i,\ell)-1}}(t+((i)\sigma_1^{(1)}-1)p^{k-1},t+p^{h(i,\ell)-1}+((i)\sigma_1^{(1)}-1)p^{k-1},\dotsc \\
&\ \ \ \ \ \ \ \ \ \ \ \ \ \ \ \ \ \ \ \ \ \ \ \ \ \ \ \ \ \ \ \ \ \ \ \ \ \ \ \ \ \ \ \ \ \ \ \ \ \ \ \ \ \ \ \ \ \ \ \ \ \ \ \ \ \ \ \ \ \dotsc,t+(p-1)p^{h(i,\ell)-1}+((i)\sigma_1^{(1)}-1)p^{k-1}) \\
&=(g_{h(i,\ell)})^{(g_k)^{(i)\sigma_1^{(1)}-1}}.
\end{align*}
This means that $\sigma_k^{(k)}$ permutes the direct factors of the base group $B$ of $P_{p^k}$ in the same way as $\sigma_1^{(1)}$ permutes the elements of $[1,p]$. In particular this implies that, for every $1\le i<k$,
\begin{equation}\label{eq:sigmakkonB}
(g_i^{(k)})^{\sigma_k^{(k)}}=(g_i^{(k-1)}, 1, \ldots, 1; 1)^{\sigma_k^{(k)}}=(1,\ldots, 1, g_i^{(k-1)}, 1, \ldots, 1; 1) \in (P_{p^{k-1}})_a \le B.
\end{equation}
Returning to \eqref{eq:star} and using \eqref{eq:sigmakkonB}, we derive:
\begin{align*}
&(f;\delta)^{\sigma_k^{(k)}}=(\prod_{i=1}^{p}\prod_{\ell=1}^{n_i}(g_{h(i,\ell)})^{(g_k)^{(i)\sigma_1^{(1)}-1}})((g_k)^{\sigma_k^{(k)}})^s \\
&=(g_{h(1,1)})^{(g_k)^{(1)\sigma_1^{(1)}-1}}\dots(g_{h(1,n_1)})^{(g_k)^{(1)\sigma_1^{(1)}-1}}\dots(g_{h(p,1)})^{(g_k)^{(p)\sigma_1^{(1)}-1}}\dots(g_{h(p,n_p)})^{(g_k)^{(p)\sigma_1^{(1)}-1}}((g_k)^a)^s \\
&=g_{h((1)\tau,1)}\dots g_{h((1)\tau,n_{(1)\tau})}\dots (g_{h((p)\tau,1)})^{(g_k)^{p-1}}\dots (g_{h((p)\tau,n_{(p)\tau})})^{(g_k)^{p-1}}((g_k)^s)^a \\
&=(\prod_{\ell=1}^{n_{(1)\tau}}g_{h((1)\tau,\ell)},1,\dotsc,1;1)\dots(1,\dotsc,1,\prod_{\ell=1}^{n_{(p)\tau}}g_{h((p)\tau,\ell)};1)(1,\dotsc,1;\gamma^s)^{a} \\
&=(f_{(1)\tau},\dotsc,f_{(p)\tau};1)(1,\dotsc,1;\delta^a)=(f_{(1)\tau},\dotsc,f_{(p)\tau};\delta^{\sigma_1^{(1)}}).
\end{align*}
The third equality holds because for distinct $i\in[1,p]$ the permutations of the form $(g_{h(i,\ell)})^{(g_k)^{(i)\sigma_1^{(1)}-1}}$ are supported on disjoint sets and thus commute. This is not the case for the last term, which acts on the full set $[1,p^k]$. The fifth equality holds by \eqref{eq:deffi}.

\noindent\textbf{Case j$<$k.} 
From \eqref{eq:star}, using (iii)-(v), we compute:
\begin{align*}
&(f_1,\dotsc,f_p;\delta)^{\sigma_j^{(k)}}=(\prod_{i=1}^{p}\prod_{\ell=1}^{n_i}(g_{h(i,\ell)})^{(g_k)^{i-1} \sigma_j^{(k)}})((g_k)^s)= (\prod_{i=1}^{p}\prod_{\ell=1}^{n_i}(g_{h(i,\ell)})^{\sigma_j^{(k)}(g_k)^{i-1}})((g_k)^s) \\
&=(g_{h(1,1)})^{\sigma_j^{(k)}}\dots(g_{h(1,n_1)})^{\sigma_j^{(k)}}\dots((g_{h(p,1)})^{\sigma_j^{(k)}})^{(g_k)^{p-1}}\dots((g_{h(p,n_p)})^{\sigma_j^{(k)}})^{(g_k)^{p-1}}((g_k)^s) \\
&=((g_{h(1,1)})^{\sigma_j^{(k-1)}},1,\dotsc,1;1)\dots((g_{h(p,n_p)})^{\sigma_j^{(k-1)}},1,\dotsc,1;1)^{(g_k)^{p-1}}(1,\dotsc,1;\gamma^s) \\
&=((\prod_{\ell=1}^{n_1}g_{h(1,\ell)})^{\sigma_j^{(k-1)}},\dotsc,(\prod_{\ell=1}^{n_p}g_{h(p,n_p)})^{\sigma_j^{(k-1)}};\gamma^s)=(f_1^{\sigma_j^{(k-1)}},\dotsc,f_p^{\sigma_j^{(k-1)}};\delta).
\end{align*}
This concludes the proof.
\end{proof}

Proposition~\ref{prop:sigmasuP} describes the action of $\langle \sigma_1^{(k)},\dotsc,\sigma_k^{(k)} \rangle$ on $P_{p^k}$ and motivates the choice of these elements as generators of the normalizer $N_{S_{p^k}}(P_{p^k})$. Indeed, for each fixed $k>1$, the $\sigma_j^{(k)}$ are elements of order $p-1$, they commute pairwise, and $(g_j^{(k)})^{\sigma_j^{(k)}} = (g_j^{(k)})^a$ for all $j\in [1,k]$. Moreover, $\sigma_j^{(k)}$ normalizes $\langle g_1^{(k)}, \dotsc, g_j^{(k)} \rangle \cong P_{p^j} \le P_{p^k}$ and centralizes $g_i^{(k)}$ for all $j<i\in[1,k]$. Accordingly, we introduce the following definition.

\begin{definition}\label{def:N_p^k}
Let $k\in\N$. We denote by $N_{p^k}$ the normalizer of $P_{p^k}$ in $S_{p^k}$, namely 
\[
N_{p^k} := N_{S_{p^k}}(P_{p^k}) = \langle g_1^{(k)},\dotsc,g_k^{(k)},\sigma_1^{(k)},\dotsc,\sigma_k^{(k)} \rangle = \langle g_i^{(k)}\mid i\in[1,k] \rangle \rtimes \langle \sigma_j^{(k)}\mid j\in[1,k]\rangle \cong P_{p^k} \rtimes (C_{p-1})^{\times k}.
\]
\end{definition}

The inherited action of $N_{p^k} / P_{p^k} \cong (C_{p-1})^{\times k}$ on $P_{p^k}$ is now fully characterized by Proposition~\ref{prop:sigmasuP}. To familiarize the reader with the notation introduced, we pause to present a concrete example, adapted from \cite{GANT}. We will return to Example~\ref{ex:111} in Section~\ref{sec:applications}, since it will be useful for the combinatorial construction of the counterexamples presented in the proofs of Propositions~\ref{rem:conclusions} and~\ref{prop:no}. 

\begin{example}\label{ex:111}
Let $p=5$, $k=2$, $n=p^k=25$ and set $a=2$. In this example we explicitly compute the elements introduced above to describe $N_{25}$. 
Here $P_{25}$ is the Sylow $5$-subgroup of $S_{25}$ generated by $g_1^{(2)}$ and $g_2^{(2)}$, where $g_1^{(2)}=g_1^{(1)}=(1,2,3,4,5)$ and 
$$g_2^{(2)}=(1,6,11,16,21)(2,7,12,17,22)(3,8,13,18,23)(4,9,14,19,24)(5,10,15,20,25).$$
Observe that $N_5=N_{S_5}(\langle g_1^{(1)}\rangle)=\langle g_1^{(1)}\rangle\rtimes\langle \sigma_1^{(1)}\rangle$, where $\sigma_1^{(1)}=(2,4,3,1)$ and $\tau=(1,3,4,2)$. Following Definitions~\ref{def:sigmajk} and~\ref{def:N_p^k}, we have that $N_{25}=P_{25}\rtimes\langle\sigma_1^{(2)}, \sigma_2^{(2)}\rangle$, where 
\begin{eqnarray*}
\sigma_1^{(2)} &=& \prod_{i=0}^4\tau_{-i5}(\sigma_1^{(1)})\tau_{i5}=(2,4,3,1)(7,9,8,6)(12,14,13,11)(17,19,18,16)(22,24,23,21), \text{and}\\
\sigma_2^{(2)} &=& \prod_{i=0}^4\tau_{i}(10,20,15,5)\tau_{-i}=(6,16,11,1)(7,17,12,2)(8,18,13,3)(9,19,14,4)(10,20,15,5).
\end{eqnarray*}
Let $\gamma$ be the $5$-cycle $(1,2,3,4,5)$ generating $P_5$, and to ease the notation let $g_1=g_1^{(2)}$ and $g_2=g_2^{(2)}$. Since $P_{25}=P_5\wr C_5=(P_5\times P_5\times P_5\times P_5\times P_5)\rtimes C_5$, then $g_1=(\gamma, 1,1,1,1;1)$, $g_2=(1,1,1,1,1;\gamma)$ and $g_1^{g_2}=(6,7,8,9,10)=(1,\gamma,1,1,1;1)$. Moreover, 
\begin{eqnarray*}
g_1^{\sigma_{1}^{(2)}} &=& (1,3,5,2,4)=(\gamma^{\sigma_{1}^{(1)}}, 1,1,1,1;1)=(\gamma^2, 1,1,1,1;1)=(g_1)^2,\\ 
g_1^{\sigma_{2}^{(2)}} &=& (6,7,8,9,10)=(1,\gamma,1,1,1;1)= g_1^{g_2},\\
g_2^{\sigma_{1}^{(2)}} &=& (1, 1,1,1,1;\gamma)=g_2,
\end{eqnarray*}
and
\begin{eqnarray*}
g_2^{\sigma_{2}^{(2)}} &=& (1,11,21,6,16)(2,12,22,7,17)(3,13,23,8,18)(4,14,24,9,19)(5,15,25,10,20)\\
&=& (1, 1,1,1,1;\gamma^2)=(g_2)^2.\\  
\end{eqnarray*}

\begin{figure}[h]
\centering
\begin{minipage}{0.6\textwidth}  
\centering
\begin{tikzpicture}[
  level distance=15mm,
  sibling distance=14mm,
  every node/.style={inner sep=1pt},
  leaf/.style={circle, fill, minimum size=3pt, inner sep=0pt}
]

\node[leaf, label=above:$\langle g_2 \rangle$] (root) {}
  child {node[leaf] (n1) {}}
  child {node[leaf] (n2) {}}
  child {node[leaf] (n3) {}}
  child {node[leaf] (n4) {}}
  child {node[leaf] (n5) {}};

\node[below=7pt of n1, anchor=north] {$\langle g_1 \rangle$};
\node[below=7pt of n2, anchor=north] {$\langle g_1^{g_2} \rangle$};
\node[below=3pt of n3, anchor=north] { \ \ $\langle g_1^{(g_2)^2} \rangle$};
\node[below=3pt of n4, anchor=north] { \ \ \ \  $\langle g_1^{(g_2)^3} \rangle$};
\node[below=3pt of n5, anchor=north] { \ \ \ \ \ \  $\langle g_1^{(g_2)^4} \rangle$};

\node at (-2,0.2) {$\sigma_2^{(2)} \longrightarrow$};
\node at (-0.3,-0.2) { \ \ \ \ \ \ \scalebox{3}{$ \hookrightarrow$}};
\node at (-4,-1.5) {$\sigma_1^{(2)} \longrightarrow$};

\end{tikzpicture}
\end{minipage}%
\hfill
\begin{minipage}{0.40\textwidth}  
\caption{The figure schematically represents the action of $\langle \sigma_1^{(2)}, \sigma_2^{(2)} \rangle$ on $P_{25}$. The element $\sigma_1^{(2)}$ centralizes $g_2$ and normalizes each $H_{i+1}=\langle g_1^{(g_2)^i} \rangle$, for $i\in [0,4]$; $\sigma_2^{(2)}$ normalizes $\langle g_2 \rangle$, and permutes the subgroups $\{ H_j \mid j\in[1,5] \}$, which are the direct factors of the base group of $P_{25}$.}
\end{minipage}
\end{figure}
\end{example}

Before proceeding, we fix some notation to be used throughout what follows.

\begin{notation}\label{not:lincp}
Given a prime $p$, consider $\omega=e^{\frac{2\pi i}{p}}$, and let $\C_p = \langle \omega\rangle$ denote the multiplicative group of complex $p$th roots of unity. Write $g=g_1^{(1)}=(1,2,\dotsc,p)$, and set $P_p=\langle g \rangle$. For each $\varepsilon\in [0,p-1]$, let $\phi_{\varepsilon}: P_p \to \C_p$ be the group homomorphism defined by $\phi_{\varepsilon}(g)=\omega^{\varepsilon}$. Then,
\begin{equation}\label{eq:lincp}
\Irr(P_p)=\Lin(P_p)=\{\phi_{\varepsilon}\mid \varepsilon\in [0,p-1] \}=\langle \phi_{1} \rangle \cong C_p.
\end{equation}
The $p$th cyclotomic field $\Q(\omega)$ contains all values of these characters and is the minimal splitting field for $P_p$. The action of $\Gal(\Q(\omega)\mid\Q)$ on $\Lin(P_p)$ is given by
\[
(\phi_{\varepsilon})^{\rho}(g)=(\phi_{\varepsilon}(g))^{\rho},	\qquad \rho\in\Gal(\Q(\omega)\mid\Q), \ \varepsilon\in [0,p-1],
\]
while the action of $N_p$ on $\Lin(P_p)$ is
\[
(\phi_{\varepsilon})^h(g)=\phi_{\varepsilon}(g^{h^{-1}})=\phi_{\varepsilon}(hgh^{-1}),	\qquad h\in N_p, \ \varepsilon\in [0,p-1].
\]
Note that $\Gal(\Q(\omega)\mid\Q)\cong\Aut(C_p)$, which is cyclic of order $p-1$, and that $\Gal(\Q(\omega)\mid\Q)$ acts transitively on $\C_p \setminus \{1\} = \{\omega^{j}\mid j\in [1,p-1] \}$ and hence also on $\Lin(P_p)\setminus\{\phi_0\}=\{\phi_{j}\mid j\in [1,p-1]\}$. In particular, these actions are naturally isomorphic. For convenience, we identify $\Gal(\Q(\omega)\mid\Q)$ with its action on the set of labels $[0,p-1]$, with $0$ as fixed point. We now explicitly associate to a generator of this group a suitable $(p-1)$-cycle acting on $[1,p-1]$.
 
Given any generator $\sigma$ of $\Gal(\Q(\omega)\mid\Q)$, let $b\in [1,p-1]$ be the primitive root modulo $p$ such that $\omega^{\sigma} = \omega^b$, and let $\pi$ be the permutation of $[0,p-1]$ such that $(\omega^{\varepsilon})^{\sigma}=(\omega^{\varepsilon})^b= \omega^{(\varepsilon){\pi}}$. Since $\phi_{\varepsilon}=(\phi_1)^{\varepsilon}$ and $(\phi_1)^{\sigma}=\phi_b$, it also follows that $(\phi_{\varepsilon})^{\sigma}=(\phi_{\varepsilon})^b=\phi_{(\varepsilon){\pi}}$ for all $\varepsilon\in [0,p-1]$. By direct computation, 
\[
\pi=(b_1,\dotsc,b_{p-1}) \qquad \text{ where } b_i\in [1,p-1] \text{ and } b_i\equiv b^i\pmod{p}, \text{ for all } i\in [1,p-1].
\]

From Definition~\ref{def:sigmajk}, let $a\in [1,p-1]$ be the primitive root modulo $p$ such that $ab \equiv 1 \pmod {p}$, so that $\sigma_1^{(1)}$ is the $(p-1)$-cycle $(a_1,\dotsc,a_{p-1})$ with $a_i\in [1,p-1]$ and $a_i\equiv a^i\pmod{p}$ for $i\in [1,p-1]$. Since $a_i \equiv a^i \equiv (b^{-1})^{i} \equiv (b^{i})^{p-1} \equiv b^{p-i} \equiv b_{p-i} \pmod{p}$ and $(\sigma_{1}^{(1)})^{-1}=(a_{p-1},\dotsc,a_1)$, it follows that $\tau=(\sigma_1^{(1)})^{-1}=(b_1,\dotsc,b_{p-1})$. Thus, $\tau$ has the same cycle decomposition as $\pi$, up to fixed points. Moreover, $\tau$ is the unique permutation of $[1,p]$ such that $g^{\tau}=(1,\dotsc,p)^{\tau}=(1,\dotsc,p)^b=g^b$ and $(p)\tau=p$. \textit{To streamline notation, we let $\tau$ act on $[0,p]$ by fixing also $0$}. In this way, for every $\varepsilon\in [0,p-1]$ one has 
\[
(\phi_{\varepsilon})^{\sigma_1^{(1)}}(g)=\phi_{\varepsilon}(g^{(\sigma_1^{(1)})^{-1}})=\phi_{\varepsilon}(g^{\tau})=\phi_{\varepsilon}(g^b)=(\phi_{\varepsilon}(g))^b=
(\phi_{\varepsilon})^b(g)=\phi_{(\varepsilon)\tau}(g)=(\phi_{\varepsilon})^{\sigma}(g), 
\]
since $(\varepsilon)\tau=(\varepsilon)\pi$. Therefore,
\begin{equation}\label{eq:t}
(\phi_{\varepsilon})^{\sigma_1^{(1)}}=\phi_{(\varepsilon)\tau}=(\phi_{\varepsilon})^{\sigma} \qquad \text{ for every } \varepsilon\in [0,p-1].
\end{equation}

In Section~\ref{sec:norgalIrr(P)} we recover this notation. Referring to $\Irr(P_p)$ as in \eqref{eq:lincp}, we proceed backwards: we first define $\sigma_1^{(1)}\in N_p$ and then choose $\sigma\in\Gal(\Q(\omega)\mid\Q)$ whose induced index permutation is $\tau = (\sigma_1^{(1)})^{-1}$. This procedure is equivalent; hence \eqref{eq:t} applies.
\end{notation} 

\subsubsection{General case}\label{sec:gensylnorgen}
What follows is based on \eqref{eq:sylowprod}.

\begin{definition}\label{def:P_n}
Let $p$ be a prime. Let $n\in \N$ and write its $p$-adic expansion as $n=\sum_{i=1}^{t}a_ip^{k_i}$, where $t\in\N$, $k_1>\cdots>k_t\geq 0$, and $a_1,\dotsc,a_t\in [1,p-1]$. We define $P_n$ the Sylow $p$-subgroup of $S_n$ given by
\[
P_n:=\prod_{i=1}^{t}\prod_{j=1}^{a_i}(P_{p^{k_i}})^{\tau(i,j)}
\]
where, for all $i\in [1,t]$, each $P_{p^{k_i}}$ denotes the Sylow $p$-subgroup of $S_{p^{k_{i}}}$ fixed in Definition~\ref{def:P_p^k}, and $\tau(i,j)\in S_n$ is defined by 
\[
\tau(i,j)=(1, 1+r(i,j))(2, 2+r(i,j))\cdots (p^{k_i}, p^{k_i}+r(i,j)),\ \text{with}\ r(i,j)=\sum_{z=1}^{i-1}a_zp^{k_z}+(j-1)p^{k_i}.
\]
\end{definition}

The following definition provides an explicit generating set for any $P_n$.

\begin{definition}
Let $p$ be a prime. Let $n\in\N$ and write $n=\sum_{i=1}^{t}a_ip^{k_i}$ as in Definition~\ref{def:P_n}. We set
\[
L(n)=\{(i,j,\ell)\in \N^{\times 3} \mid i\in [1,t],\ j\in [1,a_i],\ \ell\in [1,k_{i}]\}.
\]
For $(i,j,\ell)\in L(n)$, we define
\[
g(i,j,\ell)=(g_{\ell}^{(k_i)})^{\tau(i,j)} \in (P_{p^{k_i}})^{\tau(i,j)}\le P_n
\]
where $g_{\ell}^{(k_i)}\in P_{p^{k_i}}$, and $\tau(i,j)\in S_n$ is as in Definition~\ref{def:P_n}. Then,
\[
P_n=\langle g(i,j,\ell) \mid (i,j,\ell)\in L(n) \rangle.
\]
\end{definition}

From \cite[Lemma 4.1]{OlssonPaper}, we know that
\[
N_{S_n}(P_n)\cong N_{p^{k_1}}\wr S_{a_1}\times\cdots\times N_{p^{k_t}}\wr S_{a_t},
\] 
where $N_{p^{k_i}}$ is as in Definition~\ref{def:N_p^k}. Thanks to this result, it is possible to explicitly provide a presentation for $N_{S_n}(P_n)$. To this end, we introduce the next definition.

\begin{definition}
Let $p$ be a prime. Let $n\in\N$ and write $n=\sum_{i=1}^{t}a_ip^{k_i}$ as in Definition~\ref{def:P_n}. For all $i\in [1,t]$, we set $\overline{\xi_{i}}=\overline{\zeta_{i}}=1$ if $a_i=1$, and otherwise we set
\[
\overline{\xi_{i}}=(1, p^{k_i}+1)(2, p^{k_i}+2)\cdots (p^{k_i}, 2p^{k_i}),\ \text{and \ } \newline 
\overline{\zeta_{i}}=\prod_{j=1}^{p^{k_i}}(j, p^{k_i}+j, 2p^{k_i}+j, \dotsc, (a_{i}-1)p^{k_i}+j).
\]
Let $r_i=\sum_{e=1}^{i-1}a_{e}p^{k_e}$ and $\tau_{r_i}=(1, r_i+1)(2, r_i+2) \cdots (a_ip^{k_i}, r_i+a_ip^{k_i})$. For every $(i,j,\ell)\in L(n)$ we let
\[
\zeta_{i}=(\overline{\zeta_{i}})^{\tau_{r_i}},\ \ \ \xi_{i}=(\overline{\xi_{i}})^{\tau_{r_i}},\ \ \text{and}\ \ \sigma(i,j,\ell)=(\sigma_{\ell}^{(k_i)})^{\tau(i,j)}\in (N_{p^{k_i}})^{\tau(i,j)},
\]
where $\sigma_{\ell}^{(k_i)}\in N_{p^{k_i}}$ and $\tau(i,j)\in S_n$ is as in Definition~\ref{def:P_n}.
\end{definition}

We now describe an explicit set of generators for $N_{S_n}(P_n)$.

\begin{definition}\label{def:N_n}
Let $p$ be a prime, $n\in\N$ and $N_n := N_{S_n}(P_n)$. If $n=\sum_{i=1}^{t}a_ip^{k_i}$ as in Definition~\ref{def:P_n}, then
\begin{align*}
N_n &=
\prod_{i=1}^t\big(\big [\prod_{j=1}^{a_i}\langle P_{p^{k_i}}, \sigma_1^{(k_i)},\ldots, \sigma_{k_i}^{(k_i)}\rangle^{\tau(i,j)}]\rtimes \langle\zeta_{i},\xi_{i}\rangle\big)\\
&=
\langle g(i,j,\ell), \sigma(i,j,\ell), \zeta_{i}, \xi_{i}\ |\ i\in [1,t],\ j\in [1,a_i],\ \ell\in [1,k_i]\rangle.
\end{align*}
\end{definition}

Note that for every $i\in [1,t]$ the group $\langle \overline{\xi_{i}}, \overline{\zeta_{i}}\rangle\cong S_{a_i}$ acts as the full symmetric group on the set $\{(P_{p^{k_i}}), (P_{p^{k_i}})^{\tau(1,2)}, \ldots, (P_{p^{k_i}})^{\tau(1,a_i)}\}$. Similarly, $\langle\zeta_{i}, \xi_{i}\rangle\cong S_{a_i}$ acts as the full symmetric group on the set $\{(P_{p^{k_i}})^{\tau(i,1)}, (P_{p^{k_i}})^{\tau(i,2)}, \ldots, (P_{p^{k_i}})^{\tau(i,a_i)}\}$.
Therefore, we can define the following action.

\begin{definition}
Let $p$ be a prime, $n\in\N$ and write $n=\sum_{i=1}^{t}a_ip^{k_i}$ as in Definition~\ref{def:P_n}. Fix $i\in [1,t]$. We define an action of $\langle\zeta_{i}, \xi_{i}\rangle$ on $[1,a_i]$ by setting $x^{g}=y$ if and only if $((P_{p^{k_i}})^{\tau(i,x)})^g = (P_{p^{k_i}})^{\tau(i,y)}$, for all $g\in \langle\zeta_{i}, \xi_{i}\rangle$ and $x,y\in [1,a_i]$. 
\end{definition}

The following result extends Proposition~\ref{prop:sigmasuP} to the general case. Its proof is a straightforward consequence of Proposition~\ref{prop:sigmasuP} and of the definitions in this section, and thus is omitted.

\begin{proposition}\label{prop:NsuPgenerale}
Let $p$ be a prime, $n\in\N$ and write $n=\sum_{i=1}^{t}a_ip^{k_i}$ as in Definition~\ref{def:P_n}. Let $x\in P_n$. Then  
\[
x=\prod_{i=1}^{t}x_i, \qquad x_i=\prod_{j=1}^{a_i}(f_{(i,j)};\delta_{(i,j)})^{\tau(i,j)}
\] 
for some $(f_{(i,j)};\delta_{(i,j)})\in P_{p^{k_i}}$, $i\in[1,t]$, $j\in[1,a_i]$. Let $(i,j,\ell)\in L(n)$ and $y,z\in[1,t]$, and consider $\sigma(i,j,\ell)\in N_n$ and $\rho_y\in \langle \xi_y, \zeta_y \rangle$. Then
\begin{align*}
(x_z)^{\sigma(i,j,\ell)} &= 
\begin{cases}
x_z, & \text{if } z \neq i, \\
\displaystyle
(f_{(i,1)};\delta_{(i,1)})^{\tau(i,1)} \cdots \left((f_{(i,j)};\delta_{(i,j)})^{\sigma_{\ell}^{(k_i)}} \right)^{\tau(i,j)} \cdots 
(f_{(i,a_i)};\delta_{(i,a_i)})^{\tau(i,a_i)}, & \text{if } z = i;
\end{cases} \\
\\
(x_z)^{\rho_y} &=
\begin{cases}
x_z, & \text{if } z \neq y, \\
\displaystyle
\prod_{n=1}^{a_y}
(f_{(y,(n)\rho_y^{-1})};\delta_{(y,(n)\rho_y^{-1})})^{\tau(y,(n)\rho_y^{-1})}, & \text{if } z = y.
\end{cases}
\end{align*}
\end{proposition}

We have obtained a complete description of the action inherited by $N_n / P_n \cong \prod_{i=1}^{t}((C_{p-1})^{\times k_i})\wr S_{a_i}$ on $P_n$ for all $n\in\N$. We now illustrate the notation introduced above through the following concrete example.

\begin{example}\label{ex222}
Let $p=5$, $n=19=3\cdot 5^1+4\cdot 5^0$ and let $a=2$. In this setting we have $\tau(1,1)=1$ and
\[
\tau(1,2)=(1,6)(2,7)(3,8)(4,9)(5,10), \  
\tau(1,3)=(1,11)(2,12)(3,13)(4,14)(5,15),
\]
\[
\tau(2,1)=(1,16), \  
\tau(2,2) = (1,17), \
\tau(2,3)=(1,18),\ \text{and}\ \tau(2,4) = (1,19).
\]
In fact, $P_{19}=P_5\times (P_5)^{\tau(1,2)}\times (P_5)^{\tau(1,3)}=\langle (1,2,3,4,5), (6,7,8,9,10), (11,12,13,14,15)\rangle$ and we have
\[
L(19)= \{(1,1,1), (1,2,1), (1,3,1)\} , \ g(1,1,1)=(1,2,3,4,5), \ \sigma(1,1,1) = (2,4,3,1),
\]
\[
g(1,2,1)=(6,7,8,9,10), \ \sigma(1,2,1) = (7,9,8,6),
\]
\[
g(1,3,1)=(11,12,13,14,15),\ \sigma(1,3,1)=(12,14,13,11).
\]
Finally, we observe that: 
$$\xi_1=(1,6)(2,7)(3,8)(4,9)(5,10), \
\zeta_1=(1,6,11)(2,7,12)(3,8,13)(4,9,14)(5,10,15),$$ $$\xi_2=(16,17), \
\zeta_2=(16,17,18,19) .$$
It is clear that $\langle\xi_1, \zeta_1\rangle\cong S_{3}$ and that $\langle\xi_2, \zeta_2\rangle\cong S_{4}$. 
We conclude that 
\begin{align*}
N_{19}&=N_{S_{19}}(P_{19})=
\langle g(1,1,1), g(1,2,1), g(1,3,1), \sigma(1,1,1), \sigma(1,2,1), \sigma(1,3,1), \xi_1, \zeta_1\rangle\times\langle\xi_2, \zeta_2\rangle\\
&\cong N_{S_5}(P_5)\wr S_3\times S_4 \cong N_{S_{15}}(P_{15})\times S_4.
\end{align*}
\end{example}

\section{Irreducible characters of $P_n$ and $\cT$-functions}\label{sec:chT-functions}
In \cite{GL3}, each irreducible character of $P_n$ is associated with an equivalence class of tuples of rooted, complete, $p$-ary labelled trees. In this section we introduce the \emph{$\cT$-functions}, which encode these combinatorial structures without relying on explicit graphical representations. We establish a bijection between the irreducible characters of $P_n$ and $\cT$-functions, and use it in Sections~\ref{sec:norgalIrr(P)} and~\ref{sec:applications} to prove our main results. The following remark summarizes the aspects of the tree-based correspondence needed in the sequel; for details, see \cite[Section 3]{GL3}.

\begin{notation}\label{not:glialberi}
Let $p$ be a prime and $k\geq 1$. A rooted, complete, $p$-ary tree of height $k-1$, labelled in $[0,p]$, is a rooted tree in which each non-leaf vertex has $p$ children, each vertex is labelled in $[0,p]$, and every path from the root to a leaf has length $k-1$. We denote by $H_p^k$ the set of all such trees. The underlying unlabelled graph is called the \emph{skeleton} (see Figures~\ref{fig:trees} and~\ref{fig:skeleton} for examples). In what follows, we refer to $\Irr(P_p)$ as in \eqref{eq:lincp} of Notation~\ref{not:lincp}. For $\varepsilon\in [0,p-1]$ we abbreviate $\cX(\phi;\phi_\varepsilon)$ by $\cX(\phi;\varepsilon)$, the meaning being clear from the context. 
\end{notation}

\begin{remark}\label{rem:gl3trees}
Lemma~\ref{lem:chP-trees} characterizes $\Irr(P_{p^k})$ for every prime $p$ and $k\geq 1$. 
This result underlies the construction in \cite[Definition~3.5 and Lemma~3.6]{GL3}, where a recursive bijection between $\Irr(P_{p^k})$ and a subset of the quotient $H_p^k/_{\sim}$ is established, $\sim$ being a natural equivalence relation introduced therein. We briefly outline this correspondence. For $\theta\in\Irr(P_{p^k})$, we recover a representative of the corresponding equivalence class modulo $\sim$. With a slight abuse of notation, we denote this procedure by $t$ and write $t(\theta)$ for a representative tree obtained in this way. A rigorous presentation, together with a detailed description of the combinatorial properties of these objects, can be found in \cite{GL3}. 
 
Suppose that $k=1$ and $\Irr(P_p)=\{\phi_{\varepsilon}\mid\varepsilon\in [0,p-1] \}$. We set $t(\phi_{\varepsilon})$ to be the single-vertex tree labelled by $\varepsilon$, thus $t(\phi_{\varepsilon})=\substack{\varepsilon\\\bullet} \ $ for all $\varepsilon\in [0,p-1]$.
Assume now that $k>1$ and let $\theta\in \Irr(P_{p^k})$. We distinguish the two cases appearing in Lemma~\ref{lem:chP-trees}. If $\theta=\cX(\varphi;\varepsilon)$ for some $\varphi\in\Irr(P_{p^{k-1}})$ and $\varepsilon\in [0,p-1]$, we define $t(\theta)$ as the tree whose root is labelled by $\varepsilon$ and whose $p$-subtrees are all equal to $t(\varphi)$. Otherwise, if $\theta=\varphi_1\times\cdots\times\varphi_p\up^{P_{p^k}}$ for some mixed $\varphi_1,\dotsc,\varphi_p\in\Irr(P_{p^{k-1}})$, we define $t(\theta)$ as the tree whose root is labelled by $p$ and whose $p$ subtrees are $t(\varphi_1),\dotsc,t(\varphi_p)$, ordered from left to right. Figure~\ref{fig:extindalb} illustrates the inductive step of this procedure. 

\begin{figure}[h!]
	\centering
		\begin{subfigure}[t]{0.48\textwidth}
		\centering
		\begin{tikzpicture}[scale=1.4, every node/.style={scale=1.2}]
			\draw (0,0.5) node(R) {$\overset{\varepsilon}{\bullet}$};
			\draw (-1.6,-1.0) node[draw, rounded corners, inner sep=2pt](1) {$\begin{matrix} t(\varphi) \\ \vdots \end{matrix}$};
			\draw (-0.4,-1.0) node[draw, rounded corners, inner sep=2pt](2) {$\begin{matrix} t(\varphi)\\ \vdots \end{matrix}$};
			\draw (0.5,-0.75) node(3) {$\ldots$};
			\draw (1.6,-1.0) node[draw, rounded corners, inner sep=2pt](4) {$\begin{matrix} t(\varphi)\\ \vdots \end{matrix}$};
			\draw (R) -- (1);
			\draw (R) -- (2);
			\draw (R) -- (4);
			
		\end{tikzpicture}
		\caption{The tree $t(\theta)$ when $\theta=\cX(\varphi;\varepsilon)$, where $\varphi\in\Irr(P_{p^{k-1}})$ and $\varepsilon\in [0,p-1]$. Note that $t(\theta)\in H_p^{k}$, since $t(\varphi)$ is a rooted, complete, $p$-ary tree of height $k-2$ with labels in $[0,p]$.}
		\label{fig:extalb}
	\end{subfigure}
	\hfill	
	\begin{subfigure}[t]{0.48\textwidth}
		\centering
		\begin{tikzpicture}[scale=1.4, every node/.style={scale=1.2}]
			\draw (0,0.5) node(R) {$\overset{p}{\bullet}$};
			\draw (-1.6,-1.0) node[draw, rounded corners, inner sep=2pt](1) {$\begin{matrix} t(\varphi_1) \\ \vdots \end{matrix}$};
			\draw (-0.4,-1.0) node[draw, rounded corners, inner sep=2pt](2) {$\begin{matrix} t(\varphi_2)\\ \vdots \end{matrix}$};
			\draw (0.5,-0.75) node(3) {$\ldots$};
			\draw (1.6,-1.0) node[draw, rounded corners, inner sep=2pt](4) {$\begin{matrix} t(\varphi_p)\\ \vdots \end{matrix}$};
			\draw (R) -- (1);
			\draw (R) -- (2);
			\draw (R) -- (4);
			
		\end{tikzpicture}
		\caption{The tree $t(\theta)$ when 
$\theta=\varphi_1\times\cdots\times\varphi_p
\up^{P_{p^k}}$, for some mixed $\varphi_1,\dotsc,\varphi_p\in\Irr(P_{p^{k-1}})$. Then $t(\theta)$ depends on the chosen expression of $\theta$, since $\theta=\varphi_{(1)\pi}\times\cdots\times\varphi_{(p)\pi}\up^{P_{p^k}}$ for any $\pi\in P_p$.}
		\label{fig:indalb}
	\end{subfigure}
	
	\caption{Inductive step of the procedure $t$. 
The boxed subtrees correspond to trees of height $k-2$ arising from the previous step of the recursion.}
	\label{fig:extindalb}
\end{figure}

A priori, this construction does not yield a unique result, since $t(\theta)$ depends on the chosen expression of $\theta$. Indeed, by Remark~\ref{rem:chlimits}, any tree obtained by cyclically permuting the $p$ subtrees at each non-leaf vertex of $t(\theta)$ also corresponds to $\theta$, and arises via $t$ from a different decomposition of $\theta$. The correspondence is therefore well defined up to the following natural equivalence relation, $\sim$. Two trees are said to be \textit{equivalent} if they differ by a finite sequence of such local cyclic re-orderings of subtrees. Associating to $\theta$ the equivalence class of $t(\theta)$ yields a well-defined injective map from $\Irr(P_{p^k})$ to $H_p^k/_{\sim}$, and hence a bijection onto its image. The representatives of the classes in the image are said to be \textit{admissible trees} (see \cite[Definition 3.8]{GL3}), since they are the ones representing characters. As a consequence of \eqref{eq:pgencar}, in the general case of arbitrary $n\in\N$, each irreducible character of $P_n$ corresponds to a tuple of such equivalence classes of trees. For a visual example of this character-tree correspondence, the reader is referred to Figure~\ref{fig:trees}.
\end{remark}

\begin{figure}[h!]
	\centering
		\begin{subfigure}[t]{0.24\textwidth}
		\centering
		\begin{tikzpicture}[scale=1.0, every node/.style={scale=0.7}]
			\draw (0,0.2) node(R) {$\overset{3}{\bullet}$};
			\draw (-1.2,-0.6) node(1) {$\overset{3}{\bullet}$};
			\draw (0,-0.6) node(2) {$\overset{3}{\bullet}$};
			\draw (1.2,-0.6) node(3) {$\overset{3}{\bullet}$};
			\draw (R) -- (1);
			\draw (R) -- (2);
			\draw (R) -- (3);
			\draw (-1.5,-1.4) node(11) {$\overset{3}{\bullet}$};
			\draw (-1.2,-1.4) node(12) {$\overset{1}{\bullet}$};
			\draw (-0.9,-1.4) node(13) {$\overset{2}{\bullet}$};
			\draw (1) -- (11);
			\draw (1) -- (12);
			\draw (1) -- (13);
			\draw (-0.3,-1.4) node(21) {$\overset{2}{\bullet}$};
			\draw (0,-1.4) node(22) {$\overset{1}{\bullet}$};
			\draw (0.3,-1.4) node(23) {$\overset{2}{\bullet}$};
			\draw (2) -- (21);
			\draw (2) -- (22);
			\draw (2) -- (23);
			\draw (0.9,-1.4) node(31) {$\overset{1}{\bullet}$};
			\draw (1.2,-1.4) node(32) {$\overset{2}{\bullet}$};
			\draw (1.5,-1.4) node(33) {$\overset{1}{\bullet}$};
			\draw (3) -- (31);
			\draw (3) -- (32);
			\draw (3) -- (33);
		\end{tikzpicture}
		\caption{Inadmissible tree: the first subtree is inadmissible, because of the leaf labelled $3$. Relabelling it with a different number, we obtain an admissible tree.}
		\label{fig:inad1}
	\end{subfigure}
	\hfill	
	\begin{subfigure}[t]{0.24\textwidth}
		\centering
		\begin{tikzpicture}[scale=1.0, every node/.style={scale=0.7}]
			\draw (0,0.2) node(R) {$\overset{2}{\bullet}$};
			\draw (-1.2,-0.6) node(1) {$\overset{0}{\bullet}$};
			\draw (0,-0.6) node(2) {$\overset{0}{\bullet}$};
			\draw (1.2,-0.6) node(3) {$\overset{0}{\bullet}$};
			\draw (R) -- (1);
			\draw (R) -- (2);
			\draw (R) -- (3);
			\draw (-1.5,-1.4) node(11) {$\overset{1}{\bullet}$};
			\draw (-1.2,-1.4) node(12) {$\overset{1}{\bullet}$};
			\draw (-0.9,-1.4) node(13) {$\overset{1}{\bullet}$};
			\draw (1) -- (11);
			\draw (1) -- (12);
			\draw (1) -- (13);
			\draw (-0.3,-1.4) node(21) {$\overset{1}{\bullet}$};
			\draw (0,-1.4) node(22) {$\overset{1}{\bullet}$};
			\draw (0.3,-1.4) node(23) {$\overset{1}{\bullet}$};
			\draw (2) -- (21);
			\draw (2) -- (22);
			\draw (2) -- (23);
			\draw (0.9,-1.4) node(31) {$\overset{1}{\bullet}$};
			\draw (1.2,-1.4) node(32) {$\overset{1}{\bullet}$};
			\draw (1.5,-1.4) node(33) {$\overset{1}{\bullet}$};
			\draw (3) -- (31);
			\draw (3) -- (32);
			\draw (3) -- (33);
		\end{tikzpicture}
		\caption{Admissible tree: its equivalence class corresponds to the linear character $\cX(\cX(\phi_1;0);2)$ of $P_{27}$. 
		This tree is the unique representative of its class.}
		\label{fig:linear}
	\end{subfigure}
	\hfill		
	\begin{subfigure}[t]{0.24\textwidth}
		\centering
		\begin{tikzpicture}[scale=1.0, every node/.style={scale=0.7}]
			\draw (0,0.2) node(R) {$\overset{3}{\bullet}$};
			\draw (-1.2,-0.6) node(1) {$\overset{3}{\bullet}$};
			\draw (0,-0.6) node(2) {$\overset{0}{\bullet}$};
			\draw (1.2,-0.6) node(3) {$\overset{0}{\bullet}$};
			\draw (R) -- (1);
			\draw (R) -- (2);
			\draw (R) -- (3);
			\draw (-1.5,-1.4) node(11) {$\overset{0}{\bullet}$};
			\draw (-1.2,-1.4) node(12) {$\overset{1}{\bullet}$};
			\draw (-0.9,-1.4) node(13) {$\overset{1}{\bullet}$};
			\draw (1) -- (11);
			\draw (1) -- (12);
			\draw (1) -- (13);
			\draw (-0.3,-1.4) node(21) {$\overset{2}{\bullet}$};
			\draw (0,-1.4) node(22) {$\overset{2}{\bullet}$};
			\draw (0.3,-1.4) node(23) {$\overset{2}{\bullet}$};
			\draw (2) -- (21);
			\draw (2) -- (22);
			\draw (2) -- (23);
			\draw (0.9,-1.4) node(31) {$\overset{1}{\bullet}$};
			\draw (1.2,-1.4) node(32) {$\overset{1}{\bullet}$};
			\draw (1.5,-1.4) node(33) {$\overset{1}{\bullet}$};
			\draw (3) -- (31);
			\draw (3) -- (32);
			\draw (3) -- (33);
		\end{tikzpicture}
		\caption{Admissible tree whose equivalence class corresponds to the character \scalebox{0.9}{$\theta=\theta_1\times\theta_2\times\theta_3\up^{P_{27}}$}, where \\ \scalebox{0.9}{$\theta_1=\phi_0\times\phi_1\times\phi_1\up^{P_9}$}, \\
		$\theta_2=\cX(\phi_2;0)$, and \\
		$\theta_3=\cX(\phi_1;0)$.}
		\label{fig:tree1}
	\end{subfigure}
	\hfill
	\begin{subfigure}[t]{0.24\textwidth}
		\centering
		\begin{tikzpicture}[scale=1.0, every node/.style={scale=0.7}]
			\draw (0,0.2) node(R) {$\overset{3}{\bullet}$};
			\draw (-1.2,-0.6) node(1) {$\overset{0}{\bullet}$};
			\draw (0,-0.6) node(2) {$\overset{3}{\bullet}$};
			\draw (1.2,-0.6) node(3) {$\overset{0}{\bullet}$};
			\draw (R) -- (1);
			\draw (R) -- (2);
			\draw (R) -- (3);
			\draw (-1.5,-1.4) node(11) {$\overset{1}{\bullet}$};
			\draw (-1.2,-1.4) node(12) {$\overset{1}{\bullet}$};
			\draw (-0.9,-1.4) node(13) {$\overset{1}{\bullet}$};
			\draw (1) -- (11);
			\draw (1) -- (12);
			\draw (1) -- (13);
			\draw (-0.3,-1.4) node(21) {$\overset{1}{\bullet}$};
			\draw (0,-1.4) node(22) {$\overset{0}{\bullet}$};
			\draw (0.3,-1.4) node(23) {$\overset{1}{\bullet}$};
			\draw (2) -- (21);
			\draw (2) -- (22);
			\draw (2) -- (23);
			\draw (0.9,-1.4) node(31) {$\overset{2}{\bullet}$};
			\draw (1.2,-1.4) node(32) {$\overset{2}{\bullet}$};
			\draw (1.5,-1.4) node(33) {$\overset{2}{\bullet}$};
			\draw (3) -- (31);
			\draw (3) -- (32);
			\draw (3) -- (33);
		\end{tikzpicture}
		\caption{Admissible tree equivalent to the one in Figure~\ref{fig:tree1}. Note that \scalebox{0.9}{$\theta=\theta_3\times\theta_1\times\theta_2\up^{P_{27}}$} and \scalebox{0.9}{$\theta_1=\phi_1\times\phi_0\times\phi_1\up^{P_9}$}. This tree also represents $\theta$.}
		\label{fig:tree2}
	\end{subfigure}
	\caption{Examples of rooted, complete $3$-ary trees of height $2$, labelled in $[0,3]$. The goal is to aid the reader in understanding the notions of admissibility and equivalence of trees mentioned in Remark~\ref{rem:gl3trees}.}
	\label{fig:trees}
\end{figure}

\begin{figure}[htbp]
    \centering
    \begin{tikzpicture}[scale=1.8, every node/.style={font=\small}]
        \node (R) at (0,0.2) {$\overset{\emptyset}{\bullet}$};

        \node (1) at (-2,-0.6) {$\overset{1}{\bullet}$};
        \node (2) at (0,-0.6) {$\overset{2}{\bullet}$};
        \node (3) at (2,-0.6) {$\overset{3}{\bullet}$};
        \draw (R) -- (1);
        \draw (R) -- (2);
        \draw (R) -- (3);

        \node (11) at (-2.6,-1.2) {$\overset{11}{\bullet}$};
        \node (12) at (-2,-1.2) {$\overset{12}{\bullet}$};
        \node (13) at (-1.4,-1.2) {$\overset{13}{\bullet}$};
        \draw (1) -- (11);
        \draw (1) -- (12);
        \draw (1) -- (13);

        \node (21) at (-0.6,-1.2) {$\overset{21}{\bullet}$};
        \node (22) at (0,-1.2) {$\overset{22}{\bullet}$};
        \node (23) at (0.6,-1.2) {$\overset{23}{\bullet}$};
        \draw (2) -- (21);
        \draw (2) -- (22);
        \draw (2) -- (23);

        \node (31) at (1.4,-1.2) {$\overset{31}{\bullet}$};
        \node (32) at (2,-1.2) {$\overset{32}{\bullet}$};
        \node (33) at (2.6,-1.2) {$\overset{33}{\bullet}$};
        \draw (3) -- (31);
        \draw (3) -- (32);
        \draw (3) -- (33);

        \node (111) at (-2.8,-1.8) {\scalebox{0.7}{$\overset{111}{\bullet}$}};
        \node (112) at (-2.6,-1.8) {\scalebox{0.7}{$\overset{112}{\bullet}$}};
        \node (113) at (-2.4,-1.8) {\scalebox{0.7}{$\overset{113}{\bullet}$}};
        \draw (11) -- (111);
        \draw (11) -- (112);
        \draw (11) -- (113);

        \node (121) at (-2.2,-1.8) {\scalebox{0.7}{$\overset{121}{\bullet}$}};
        \node (122) at (-2,-1.8) {\scalebox{0.7}{$\overset{122}{\bullet}$}};
        \node (123) at (-1.8,-1.8) {\scalebox{0.7}{$\overset{123}{\bullet}$}};
        \draw (12) -- (121);
        \draw (12) -- (122);
        \draw (12) -- (123);

        \node (131) at (-1.6,-1.8) {\scalebox{0.7}{$\overset{131}{\bullet}$}};
        \node (132) at (-1.4,-1.8) {\scalebox{0.7}{$\overset{132}{\bullet}$}};
        \node (133) at (-1.2,-1.8) {\scalebox{0.7}{$\overset{133}{\bullet}$}};
        \draw (13) -- (131);
        \draw (13) -- (132);
        \draw (13) -- (133);

        \node (211) at (-0.8,-1.8) {\scalebox{0.7}{$\overset{211}{\bullet}$}};
        \node (212) at (-0.6,-1.8) {\scalebox{0.7}{$\overset{212}{\bullet}$}};
        \node (213) at (-0.4,-1.8) {\scalebox{0.7}{$\overset{213}{\bullet}$}};
        \draw (21) -- (211);
        \draw (21) -- (212);
        \draw (21) -- (213);

        \node (221) at (-0.2,-1.8) {\scalebox{0.7}{$\overset{221}{\bullet}$}};
        \node (222) at (0,-1.8) {\scalebox{0.7}{$\overset{222}{\bullet}$}};
        \node (223) at (0.2,-1.8) {\scalebox{0.7}{$\overset{223}{\bullet}$}};
        \draw (22) -- (221);
        \draw (22) -- (222);
        \draw (22) -- (223);

        \node (231) at (0.4,-1.8) {\scalebox{0.7}{$\overset{231}{\bullet}$}};
        \node (232) at (0.6,-1.8) {\scalebox{0.7}{$\overset{232}{\bullet}$}};
        \node (233) at (0.8,-1.8) {\scalebox{0.7}{$\overset{233}{\bullet}$}};
        \draw (23) -- (231);
        \draw (23) -- (232);
        \draw (23) -- (233);

        \node (311) at (1.2,-1.8) {\scalebox{0.7}{$\overset{311}{\bullet}$}};
        \node (312) at (1.4,-1.8) {\scalebox{0.7}{$\overset{312}{\bullet}$}};
        \node (313) at (1.6,-1.8) {\scalebox{0.7}{$\overset{313}{\bullet}$}};
        \draw (31) -- (311);
        \draw (31) -- (312);
        \draw (31) -- (313);

        \node (321) at (1.8,-1.8) {\scalebox{0.7}{$\overset{321}{\bullet}$}};
        \node (322) at (2,-1.8) {\scalebox{0.7}{$\overset{322}{\bullet}$}};
        \node (323) at (2.2,-1.8) {\scalebox{0.7}{$\overset{323}{\bullet}$}};
        \draw (32) -- (321);
        \draw (32) -- (322);
        \draw (32) -- (323);

        \node (331) at (2.4,-1.8) {\scalebox{0.7}{$\overset{331}{\bullet}$}};
        \node (332) at (2.6,-1.8) {\scalebox{0.7}{$\overset{332}{\bullet}$}};
        \node (333) at (2.8,-1.8) {\scalebox{0.7}{$\overset{333}{\bullet}$}};
        \draw (33) -- (331);
        \draw (33) -- (332);
        \draw (33) -- (333);
    \end{tikzpicture}
    \caption{Visualization of \scalebox{1.2}{$\sk_{3^4}$}, representing the skeleton of the rooted, complete, $3$-ary tree of \mbox{height $3$}. The unnamed version of the graph depicted here is the skeleton referenced in Remark~\ref{rem:gl3trees}.
    }
    \label{fig:skeleton}
\end{figure}

We are now ready to proceed with the introduction of the new formalism. 

\begin{notation}\label{not:funcsequences}
For any two sets $A$ and $B$, it is standard to denote by $B^A$ the set of all functions from $A$ to $B$. In the specific case where $A=[1,\ell]$ for some $\ell\in\N$, we may equivalently regard $B^{[1,\ell]}$ as the set of sequences of length $\ell$ with entries in $B$. Note that $[1,\ell]=\emptyset$ when $\ell=0$, and that the only element of $B^{\emptyset}$ is the empty sequence, also denoted by $\emptyset$; if $B=\emptyset$ and $\ell > 0$ instead, then $B^{[1,\ell]}$ is the empty set. Given $\ell > 0$ and any $\ts\in B^{[1,\ell]}$ we write $s_i$ for $\ts(i)$, for all $i\in[1,\ell]$.
\end{notation}

\begin{definition}\label{def:sk} 
Given a prime $p$ and an integer $k\geq 1$, $$\sk_{p^k}:=\bigcup_{i=0}^{k-1}[1,p]^{[1,i]}$$ is the disjoint union of sequences of any length $i\in[0,k-1]$ with entries in $[1,p]$. For any $\ell\in[1,k-1]$, a sequence $\ts\in\sk_{p^k}$ of length $\ell(\ts)=\ell$ is denoted as $\ts=s_1\dots s_{\ell}$, where each $s_j\in[1,p]$ for all $j\in[1,\ell]$. Given two sequences $\mathtt{g},\mathtt{h}\in\sk_{p^k}$, if $\ell(\mathtt{g})+\ell(\mathtt{h})\le k-1$ the sequence $\mathtt{g} |\mathtt{h}:=\mathtt{gh}$ obtained by concatenation (that is, appending $\mathtt{h}$ to $\mathtt{g}$) is also a sequence of $\sk_{p^k}$. Moreover, $\mathtt{h}$ is a child of $\mathtt{g}$, in symbols $\mathtt{g}\sqsubseteq \mathtt{h}$, if $\mathtt{h}=\mathtt{g}x$ for some $x\in[1,p]$. This is a partial order on $\sk_{p^k}$, and $(\sk_{p^k},\sqsubseteq)$ is called the $p^k$-skeleton.
\end{definition}

Figure~\ref{fig:skeleton} shows how ordered sets of sequences can be used to represent graphs of unlabelled trees, thereby motivating the introduction of the $p^k$-skeleton to describe the skeleton of a rooted, complete, $p$-ary tree of height $k-1$ (as in Notation~\ref{not:glialberi}). Observe that if $k=1$ then $\sk_p=\{\emptyset\}$ contains only the empty sequence.

\begin{notation}\label{not:denotingsequences}
Let $j\in\N$, and let $\ts\in\sk_{p^k}$ be a sequence of length $\ell\in[1,k-1]$. We denote by $\ts^{j}$ the prefix of $\ts$ consisting of its first $j$ symbols, namely $\ts^{j}:=s_1\dots s_j$. We denote by $^{j}\ts$ the suffix of $\ts$ obtained by removing its first $j$ symbols, namely $^{j}\ts:=s_{j+1}\dots s_{\ell}$. Accordingly, we set $\ts^{0}=\emptyset$, $^{0}\ts=\ts$ and $\ts^{j}=\ts$, $^{j}\ts=\emptyset$ if $j>\ell$. In particular, $\ts=\ts^{j}|^{j}\ts$ for every $j\in\N$.
\end{notation}

Rooted complete $p$-ary trees of height $k-1$ labelled in $[0,p]$ are identified with functions from the $p^k$-skeleton to $[0,p]$. 

\begin{definition}\label{def:labelingfunctions}
Let $p$ be a prime and $k\geq 1$ an integer. We denote by $\tF_{p^k}$ the set defined by $$\tF_{p^k}:=[0,p]^{\sk_{p^k}}.$$ A generic element $t : \sk_{p^k} \to [0,p]$ of this set is called a $p^k$-\textit{labeling} function.
\end{definition}

\begin{example}\label{ex:datreealabelfunc}
We consider the $3^3$-labeling functions associated with the trees in Figures~\ref{fig:trees} and denote them by $t_a, t_b, t_c, t_d$, respectively.
\begin{itemize}
	\item[(A)] $t_a^{-1}(\{0\})=\emptyset,	\ t_a^{-1}(\{1\})=\{ 12, 22, 31, 33 \}, \	t_a^{-1}(\{2\})=\{ 13, 21, 23, 32 \}, \ 	t_a^{-1}(\{3\})=\{ \emptyset, 1, 2, 3, 11 \}$.
	\item[(B)] $t_b^{-1}(\{0\})=\{ 1, 2, 3 \},	\ t_b^{-1}(\{1\})=\{ 11, 12, 13, 21, 22, 23, 31, 32, 33 \},	\ t_b^{-1}(\{2\})=\{ \emptyset \}, \	t_b^{-1}(\{3\})=\emptyset$. 
	\item[(C)] $t_c^{-1}(\{0\})=\{ 2, 3, 11 \},	\ t_c^{-1}(\{1\})=\{ 12, 13, 31, 32, 33 \}, \ 	t_c^{-1}(\{2\})=\{ 21, 22, 23 \},	\ t_c^{-1}(\{3\})=\{ \emptyset, 1 \}$.
	\item[(D)] $t_d^{-1}(\{0\})=\{ 1, 3, 22 \},	\ t_d^{-1}(\{1\})=\{ 11, 12, 13, 21, 23 \}, \	t_d^{-1}(\{2\})=\{ 31, 32, 33 \}, \ 	t_d^{-1}(\{3\})=\{ \emptyset, 2 \}$.
\end{itemize} 
\end{example}

The following definition encodes the equivalence relation on rooted, complete, $p$-ary labelled trees (see Remark~\ref{rem:gl3trees}) in terms of $p^k$-labelling functions. 

\begin{definition}\label{def:ppdctreefunc}
Let $p$ be a prime and $k\geq 1$ an integer. For $t_1,t_2\in\tF_{p^k}$, we write $t_1 \sim_{p^k} t_2$ if and only if $t_1(\emptyset)=t_2(\emptyset)$ and for every $\mathtt{q}\in\sk_{p^{k-1}}$ there exist $\gamma_{\mathtt{q}}\in P_p$ such that 
\[
t_1(\ts)=t_2((s_1)\gamma_{\ts^0} \dots (s_i)\gamma_{\ts^{i-1}} \dots (s_{\ell(\ts)})\gamma_{\ts^{\ell(\ts)-1}}) \qquad \text{ for any } \ts\in\sk_{p^k} \smallsetminus \{ \emptyset \}.
\]
This is an equivalence relation $\sim_{p^k}$ on $\tF_{p^k}$. We denote by $\bF_{p^k}:= \tF_{p^k}/\!\sim_{p^k}$ the corresponding quotient set. We refer to the elements of $\bF_{p^k}$ as $p^k$-\textit{tree functions}.
\end{definition}

\begin{notation}\label{not:representatives}
To simplify notation, we identify each element of $\bF_{p^k}$ with a representative $h\in\tF_{p^k}$ and denote it by the corresponding capital letter $H$. For any $\ts\in\sk_{p^k}$, we set $H(\ts)=h(\ts)$. All results are independent of the choice of representatives.
\end{notation}

\begin{example}\label{ex:fromttoT}
For a concrete display of Definition~\ref{def:ppdctreefunc} and Notation~\ref{not:representatives} we return to Figure~\ref{fig:trees} and Example~\ref{ex:datreealabelfunc}. We prove that $t_c$ and $t_d$ belong to the same equivalence class. Let $\gamma_{\emptyset}=\gamma_{1}=(1,2,3)$ and $\gamma_{2}=\gamma_{3}=1_{P_3}$. Then $t_c(\emptyset)=t_d(\emptyset)$, and
\begin{align*}
&t_c(1)=t_d((1)\gamma_{\emptyset})=t_d(2)=3, \
t_c(2)=t_d((2)\gamma_{\emptyset})=t_d(3)=0, \
t_c(3)=t_d((3)\gamma_{\emptyset})=t_d(1)=0, \\
&t_c(11)=t_d((1)\gamma_{\emptyset}\mid(1)\gamma_{1})=t_d(22)=0, \
t_c(12)=t_d((1)\gamma_{\emptyset}\mid(2)\gamma_{1})=t_d(23)=1, \\
&t_c(13)=t_d((1)\gamma_{\emptyset}\mid(3)\gamma_{1})=t_d(21)=1, \
t_c(21)=t_d((2)\gamma_{\emptyset}\mid(1)\gamma_{2})=t_d(31)=2, \\
&t_c(22)=t_d((2)\gamma_{\emptyset}\mid(2)\gamma_{2})=t_d(32)=2, \
t_c(23)=t_d((2)\gamma_{\emptyset}\mid(3)\gamma_{2})=t_d(33)=2, \\
&t_c(31)=t_d((3)\gamma_{\emptyset}\mid(1)\gamma_{3})=t_d(11)=1, \
t_c(32)=t_d((3)\gamma_{\emptyset}\mid(2)\gamma_{3})=t_d(12)=1, \\
&t_c(33)=t_d((3)\gamma_{\emptyset}\mid(3)\gamma_{3})=t_d(13)=1.
\end{align*}
Therefore, $t_c\sim_{3^3}t_d$.  Observe that the tree corresponding to $t_d$ is obtained from that of $t_c$ by cyclically permuting the three subtrees below the root, and then the leaves of the height one subtree whose root has label $3$. This reasoning suggests that $T_b=\{t_b\}$. In fact, suppose that $t\in F_{3^3}$ is such that $t\sim_{3^3}t_b$. As a direct consequence of Definition~\ref{def:ppdctreefunc}, it follows that $t(\emptyset)=2$ and $t(\mathtt{g})=0$, $t(\mathtt{h})=1$ for any $\mathtt{g}, \mathtt{h}\in\sk_{3^3}$ with $\ell(\mathtt{g})=1$ and $\ell(\mathtt{h})=2$. This implies $t=t_b$. If $t_c$ is chosen to represent the class of $t_c$ (and $t_d$), we may denote such class by $T_c$ and adopt shorthand expressions like $T_c(22)=t_c(22)=2$.
\end{example}

We now consider subtrees rooted at arbitrary vertices.

\begin{definition}\label{def:subfunctions}
Let $p$ be a prime, $k\geq 1$ an integer and $t\in\tF_{p^k}$. For any $\ts\in\sk_{p^k}$, define $(t(\ts))^{\giu}\in\tF_{p^{k- \ell(\ts)}}$ by 
\[ 
(t(\ts))^{\giu}(\mathtt{q})=t(\mathtt{sq}), \ \text{ \ for all \ } \mathtt{q}\in\sk_{p^{k- \ell(\ts)}}.
\] 
We call $(t(\ts))^{\giu}$ the $p^{k- \ell(\ts)}$-\textit{labeling subfunction of $t$ rooted at $\ts$}. We use the symbol $(T(\ts))^{\giu}$ to denote its equivalence class in $\bF_{p^{k- \ell(s)}}$, referred to as the \textit{tree subfunction}  of $t$ rooted at $\ts$. 
\end{definition}

\begin{example}\label{ex:subfunc}
Consider again Figure~\ref{fig:trees} and Example~\ref{ex:datreealabelfunc}. Note that $\sk_{3^2}=\{\emptyset,1,2,3\}$. We compute the $3^2$-subfunction of $t_c$ rooted at $1$, $(t_c(1))^{\giu}$, and the $3^2$-subfunction of $t_d$ rooted at $2$, $(t_d(2))^{\giu}$. Then, 
\begin{align*}
&(t_c(1))^{\giu}(\emptyset)=t_c(1)=3, \
(t_c(1))^{\giu}(1)=t_c(11)=0, \
(t_c(1))^{\giu}(2)=t_c(12)=1, \
(t_c(1))^{\giu}(3)=t_c(13)=1, \\
&(t_d(2))^{\giu}(\emptyset)=t_d(2)=3, \
(t_d(2))^{\giu}(1)=t_d(21)=1, \
(t_d(2))^{\giu}(2)=t_d(22)=0, \
(t_d(2))^{\giu}(3)=t_d(23)=1.
\end{align*}
Observe that $(t_c(1))^{\giu},(t_d(2))^{\giu}\in \tF_{3^2}$ and that $(t_c(1))^{\giu}\sim_{3^2}(t_d(2))^{\giu}$. Hence, by Definition~\ref{def:subfunctions} and Notation~\ref{not:representatives}, we may take $(t_c(1))^{\giu}$ as a representative of its class $(T_c(1))^{\giu}\in\bF_{3^2}$ and write, for instance, $(T_c(1))^{\giu}(2)=(t_c(1))^{\giu}(2)=0$. Note that $(t_c(\emptyset))^{\giu}=t_c$ and $(T_c(\emptyset))^{\giu}=T_c$.
\end{example}

\begin{notation}\label{not:T=()}
Let $k\geq 2$. For $T\in\bF_{p^k}$ and $\ts\in\sk_{p^{k-1}}$, the elements $(T(\ts 1))^{\giu},\dotsc,(T(\ts p))^{\giu}\in\bF_{p^{k- (\ell(s)+1)}}$ are the \textit{tree subfunctions} of $(T(\ts))^{\giu}$ rooted at the children of $\ts$. Observe that $T$ is uniquely determined by $T(\emptyset)$ and by its $p^{k-1}$-tree subfunctions $(T(1))^{\giu},\dotsc,(T(p))^{\giu}$, up to a cyclic permutation of $P_p$. Accordingly, setting $x=T(\emptyset)$ and $T_i=(T(i))^{\giu}$ for every $i\in [1,p]$, we represent $T$ by 
\[
T=(T_1\mid\dots\mid T_p;x).
\]  
Conversely, for every $i\in[1,p]$ consider $f_i\in\tF_{p^{k-1}}$ and the corresponding $F_i\in\bF_{p^{k-1}}$. Let $\varepsilon\in[0,p]$ and write $f=(f_1 \mid \dots \mid f_p;\varepsilon)$ for the element of $\tF_{p^k}$ such that $f(\emptyset)=\varepsilon$, and $(f(i))^{\giu}=f_i$ for each $i\in[1,p]$. We denote by \[
F=(F_1 \mid \dots \mid F_p;\varepsilon)
\] 
the $p^k$-tree function associated to $f$. If $k=1$, these expressions simply reduce to $T=(x)$ and $F=(\varepsilon)$.
\end{notation}

\begin{example}
We refer to Figure~\ref{fig:tree1} and Examples~\ref{ex:datreealabelfunc} and \ref{ex:subfunc} to give an example of Notation~\ref{not:T=()}. Consider $(t_c(i))^{\giu}\in\tF_{3^2}$ and $(T_c(i))^{\giu}\in\bF_{3^2}$ for all $i\in[1,3]$. We write $t_c=((t_c(1))^{\giu}\mid(t_c(2))^{\giu}\mid(t_c(3))^{\giu};3)$ and $T_c=((T_c(1))^{\giu}\mid(T_c(2))^{\giu}\mid(T_c(3))^{\giu};3)$.
\end{example}

We now focus on a distinguished subset of $\bF_{p^k}$. The notion of tree admissibility (see \cite[Definition 3.8]{GL3}) is reformulated in terms of $\cT$-functions. 

\begin{definition}\label{def:admT}
Let $p$ be a prime, $k\geq 1$ an integer and $T\in\bF_{p^k}$. We call $T$ an \textit{admissible $p^k$-tree function} if it satisfies the following properties. If $k=1$, $T$ is admissible if and only if $T(\emptyset)\neq p$. If $k>1$, $T$ is admissible if and only if for any $\ts\in\sk_{p^{k-1}}$ exactly one of the following holds:
\begin{itemize}
	\item[(i)] $T(\ts)\neq p$ and $(T(\ts 1))^{\giu},\dotsc, (T(\ts p))^{\giu}$ are admissible and equal; 
	\item[(ii)] $T(\ts)=p$ and $(T(\ts 1))^{\giu},\dotsc,(T(\ts p))^{\giu}$ are admissible and mixed. 
\end{itemize}
If $T$ is admissible, we call it a $\cT$-function. We denote by $\cF_{p^k}$ the specific subset of $\bF_{p^k}$ consisting of all $\cT$-functions.
\end{definition}

\begin{notation}\label{not:cTfunc}
We say that a function $t\in\tF_{p^k}$ is an \emph{admissible} $p^k$-labeling function if it represents an admissible $p^k$-tree function in $\cF_{p^k}$. In this case, we denote its equivalence class by the corresponding calligraphic letter $\cT$. This convention applies only to admissible functions. Accordingly, elements of $\cF_{p^k}$ are called $\cT$-functions. For $k=1$, we have $\tF_p=\{t^j\mid j\in[0,p]\}$, where $t^j:\sk_p\to [0,p]$ is uniquely determined by $t^j(\emptyset)=j$. The admissible $p$-labeling functions are $t^0,\dotsc,t^{p-1}$, while $t^p$ is not admissible since $t^p(\emptyset)=p$. Hence, $\cF_p=\{\cT^{\varepsilon}\mid \varepsilon\in [0,p-1]\}$ where $\cT^{\varepsilon}=(\varepsilon)$ by Notation~\ref{not:T=()}.
\end{notation}

\begin{example}
We return to Figure~\ref{fig:trees} and Examples~\ref{ex:datreealabelfunc} and \ref{ex:fromttoT}. Note that $(T_a(1))^{\giu}$ is an inadmissible $3^2$-tree function, since $(T_a(11))^{\giu}$ is not admissible. Hence $T_a$ is not a $\cT$-function. Equivalently, the tree in Figure~\ref{fig:inad1} is not admissible, as its subtree rooted at the first vertex below the root is inadmissible according to Remark~\ref{rem:gl3trees}. The remaining $3^3$-tree functions $T_b,T_c$ are admissible. By Notation~\ref{not:cTfunc}, we denote them by $\cT_B$ and $\cT_C$, respectively, thereby implicitly conveying their admissibility. This convention will be used throughout without further comment. 
\end{example}

Equivalence classes of admissible trees and $\cT$-functions provide two equivalent descriptions of the same objects, and are in natural bijection via the identification introduced in this section. We thus freely identify trees and functions throughout. Using Definition~\ref{def:admT}, we reformulate the character-tree correspondence of \cite{GL3} as follows, obtaining a bijection between irreducible characters and $\cT$-functions. 

\begin{definition}(\cite[Definitions 3.5 and 3.8]{GL3})\label{def:kcorr}
Let $p$ be a prime and $k\geq 1$ an integer.
We define recursively a map 
		\[ \Phi_{p^k}:\Irr(P_{p^k})\to\cF_{p^k} \]
		as follows:
		\begin{itemize}
			\item[(i)] for $k=1$ and any $\varepsilon\in[0,p-1]$, set $\Phi_p(\phi_\varepsilon)=\cT^{\varepsilon}$;
			
			\item[(ii)] for any $k\ge 2$ and $\theta\in\Irr(P_{p^k})$, set
			\[ \Phi_{p^k}(\theta) = \begin{cases}
				(\cT_1 \mid \dots \mid \cT_p; \varepsilon) & \text{if } \theta=\cX(\varphi;\varepsilon)\text{ for some }\varphi\in\Irr(P_{p^{k-1}})\text{ and }\varepsilon\in[0,p-1],\\
				(\cT_1 \mid \dots \mid \cT_p; p)  & \text{if } \theta=\varphi_1\times\cdots\times\varphi_p\up^{P_{p^k}} \text{ for some mixed }\varphi_1,\dotsc,\varphi_p\in\Irr(P_{p^{k-1}}),
			\end{cases} \]
			where $\cT_i=\Phi_{p^{k-1}}(\varphi)$ in the first case, and $\cT_i=\Phi_{p^{k-1}}(\varphi_i)$ in the second case, for all $i\in[1,p]$. 
		\end{itemize}
\end{definition}

This construction closely mirrors the procedure in Remark~\ref{rem:gl3trees} and agrees with the correspondence in \cite{GL3}. We now describe the inverse map.

\begin{lemma}\textnormal{(\cite[Lemma 3.6]{GL3})}\label{lem:Tbij}
Let $p$ be a prime and $k\geq 1$ an integer. The map $\Phi_{p^k}$ is well defined and a bijection between $\Irr(P_{p^k})$ and $\cF_{p^k}$. Its inverse $\Psi_{p^k}: \cF_{p^k} \longrightarrow \Irr(P_{p^k})$ is uniquely determined by 
\[
\Psi_{p^k}\circ\Phi_{p^k}=\mathrm{Id}_{\Irr(P_{p^k})}, \qquad \Phi_{p^k}\circ\Psi_{p^k}=\mathrm{Id}_{\cF_{p^k}}
\]
and is defined recursively as follows:
\begin{itemize}
	\item[(i)] if $k=1$, then $\Psi_p(\cT^{\varepsilon})=\phi_\varepsilon$ for every $\varepsilon\in[0,p-1]$;
	
	\item[(ii)] for any $k\geq 2$ and any $\cT=(\cT_1\mid\dots\mid \cT_p;\varepsilon)\in\cF_{p^k}$, one has
\[ \Psi_{p^k}(\cT) = \begin{cases}
	\cX\big(\Psi_{p^{k-1}}(\cT_1); \varepsilon\big) & \text{if }\varepsilon\in[0,p-1],\\
	\Psi_{p^{k-1}}(\cT_1) \times \cdots \times \Psi_{p^{k-1}}(\cT_p) \up^{P_{p^k}} & \text{if }\varepsilon=p.
\end{cases} \]
\end{itemize}
We write $\theta(\cT):=\Psi_{p^k}(\cT)$ for the irreducible character of $P_{p^k}$ corresponding to $\cT$. Similarly, given $\theta\in\Irr(P_{p^k})$, we write $\cT(\theta):=\Phi_{p^k}(\theta)$ for the associated $\cT$-function.
\end{lemma}

This formalism provides an equivalent framework for our study. We now extend it to the general case. 

\begin{notation}\label{not:p-adicn}
Let $p$ be a prime. Given any $n\in\N$, we write its $p$-adic expansion as 
\[
n=\sum_{i=1}^t a_ip^{k_i} \qquad \text{for some } t\in\N, \text{ with integers } k_1>\ldots>k_t\geq 0 \text{ and } a_1,\dotsc,a_t\in [1,p-1].
\]
Set $q_0=0$ and $q_{i+1}=q_i+a_{i+1}$ for all $i\in[0,t-1]$. Let $\mathtt{o}: [1,q_t] \to \{k_1,\dotsc,k_t\}$ be the function defined by 
\[
\mathtt{o}(j)=k_i \qquad  \text{ if } q_{i-1}< j \le q_i.
\]
The decompositions
\begin{align*}
P_n\cong (P_{p^{k_1}})^{\times a_1}\times\cdots\times (P_{p^{k_t}})^{\times a_t}, \qquad
\Irr(P_n)= \{\theta_1\times\cdots\times\theta_{q_t} \mid  \theta_j\in\Irr(P_{p^{\mathtt{o}(j)}}) \ \text{for } j\in [1,q_t]\}
\end{align*}
are understood up to trivial factors, which may be disregarded. Indeed, if $k_t=0$ then $(P_{p^{k_t}})^{\times a_t}$ is trivial and contributes only the trivial character. In what follows, the terms corresponding to $k_t$ may be safely ignored when $k_t=0$.  
\end{notation}

We now generalize Definitions~\ref{def:sk}--\ref{def:kcorr} in the sense of \cite[Definition 3.19]{GL3}.

\begin{definition}\label{def:cTfuncgen}
Let $p$ be a prime and $n\in\N$ with $p$-adic expansion as in Notation~\ref{not:p-adicn}. The $n$-skeleton is 
\[
\sk_n:=(\sk_{p^{k_1}})^{\times a_1} \times \dots \times (\sk_{p^{k_t}})^{\times a_t}
\] 
and 
\[
\tF_n:=(\tF_{p^{k_1}})^{\times a_1} \times \cdots \times (\tF_{p^{k_t}})^{\times a_t}
\]
is the set of $n$-labeling functions. As a result, $\mathbf{t}\in\tF_n$ is a tuple of $t_j\in\tF_{p^{\mathtt{o}(j)}}$ for all $j\in[1,q_t]$, so that $\mathbf{t}=(t_1,\dotsc,t_{q_t})$ can be viewed as a function of several variables. Given two elements $\mathbf{x},\mathbf{y}\in \tF_n$, we set $\mathbf{x}\sim_n \mathbf{y}$ if and only if $X_j=Y_j$, for all $j\in[1,q_t]$. The quotient set is given by
\[
\bF_n=(\bF_{p^{k_1}})^{\times a_1} \times \cdots \times (\bF_{p^{k_t}})^{ \times a_t},
\]
and the set of $\cT$-functions is naturally defined as 
\[
\cF_n=(\cF_{p^{k_1}})^{\times a_1} \times \dots \times (\cF_{p^{k_t}})^{\times a_t}.
\]
With no risk of ambiguity, we shall refer to $\cT\in\cF_n$ as an admissible $n$-tree function, and denote it by 
\[
\cT=(\cT_1,\dotsc,\cT_{q_t}), \qquad \text{with } \ \cT_j\in\cF_{p^{\mathtt{o}(j)}} \ \text{ for all } j\in[1,q_t].
\] 
\end{definition}

We freely identify $\cT=(\cT_1,\dotsc,\cT_{q_t})\in\cF_n$ with a representative $\mathbf{t}=(t_1,\dotsc,t_{q_t})\in\tF_n$. Given $\mathbf{\ts}\in\sk_n$, we denote by $\cT(\mathbf{\ts})=(\cT_1(\ts_1),\dotsc,\cT_{q_t}(\ts_{q_t}))$ the element $\mathbf{t}(\mathbf{\ts})=(t_1(\ts_1),\dotsc,t_{q_t}(\ts_{q_t}))\in [0,p]^{\times q_t}$, where $\ts_{j}$ is the $j$th entry of the tuple $\mathbf{\ts}$. For example, for any $n\in\N$, we have $\cT(\triv_{P_n})(\mathbf{\ts})=(0,\dotsc,0)$ for all $\mathbf{\ts}\in\sk_n$.

\begin{definition}\label{def:n-corr} \textnormal{(\cite[Definition 3.19]{GL3})}
Let $p$ be a prime and $n\in\N$ with $p$-adic expansion as in Notation~\ref{not:p-adicn}. We define the bijection $\Phi_n : \Irr(P_n) \to \cF_n$ by
\[
\Phi_n(\theta):=(\Phi_{p^{\mathtt{o}(1)}}(\theta_1),\dotsc,\Phi_{p^{\mathtt{o}(j)}}(\theta_j),\dotsc,\Phi_{p^{\mathtt{o}(q_t)}}(\theta_{q_t})), \qquad \text{ for all } \theta\in\Irr(P_n).
\]
Similarly, the inverse map $\Psi_n : \cF_n \to \Irr(P_n)$ is defined by
\[
\Psi_n(\cT)=\Psi_{p^{\mathtt{o}(1)}}(\cT_1)\times\cdots\times\Psi_{p^{\mathtt{o}(j)}}(\cT_j)\times\cdots\times\Psi_{p^{\mathtt{o}(q_t)}}(\cT_{q_t}), \qquad \text{ for all } \cT\in\cF_n.
\]
Then, with a mild abuse of notation, $\cT(\theta)=(\cT(\theta_1),\dotsc,\cT(\theta_{q_t}))$ denotes the $\cT$-function corresponding to $\theta$, while
$\theta(\cT)=\theta(\cT_1)\times\cdots\times\theta(\cT_{q_t})$ denotes the irreducible character of $P_n$ corresponding to $\cT$.
\end{definition}

This completes the translation of the character-tree correspondence of \cite{GL3} in terms of \mbox{$\cT$-functions}, which will now serve as the sole language for the analysis that follows. In conclusion, \textit{for any prime $p$ and $n\in\N$, each irreducible character $\theta$ of $P_n$ is uniquely determined and described by $\cT(\theta)$.}

The \textit{admissible tree-statistics} (or \textit{character-statistics}) from \cite[Section 3]{GL3}, admit a direct reformulation as \textit{$\cT$-statistics}, along with the associated results. These statistics are used to describe the positivity of Sylow branching coefficients for symmetric groups, namely the multiplicities $[\chi\down_{P_n},\phi]$ where $\chi$ and $\phi$ are any irreducible characters of $S_n$ and $P_n$, respectively. We do not develop this further here, but we include a reformulation of \cite[Proposition 3.20]{GL3} in terms of $\cT$-functions.

\begin{definition}\label{def:retrotreefunc}
Let $p$ be a prime and $k\geq 1$ an integer. 
Given $\cT\in\cF_{p^k}$ and $i\in[0,p]$, we define $(\cT)^{-1}(i):=\{ \ts\in\sk_{p^k} \mid \cT(\ts)=i \}$.
\end{definition}

The sets above partition the entries (of a representative) of a $\cT$-function by label. The number of labels $p$ reflects a specific algebraic property of the corresponding irreducible character.

\begin{proposition}\label{prop:degree}
\textnormal{(\cite[Proposition 3.20]{GL3})}
Let $p$ be a prime and $n\in\N$ with $p$-adic expansion as in Notation~\ref{not:p-adicn}. Let $\theta\in\Irr(P_n)$ and $\cT=\cT(\theta)=(\cT_1,\dotsc,\cT_{q_t})$. If $m=\sum_{i=1}^{q_t} \left| (\cT_i)^{-1}(p) \right|$, then $\theta(1)=p^m$.
\end{proposition}

\begin{example}
The irreducible character of Figure~\ref{fig:tree1} has degree $3^2=9$, while the irreducible character of Figure~\ref{fig:indno0p} has degree $5^2=25$.
\end{example}

\section{Applications}\label{sec:norgalIrr(P)}
Throughout this section, we focus on the fixed Sylow $p$-subgroup $P_n$ and its normalizer $N_n$ described in Subsection~\ref{sec:gensylnor}. We refer to irreducible characters of $P_n$ via $\cT$-functions, as explained in Section~\ref{sec:chT-functions}. More precisely, we use the parametrization given in Definitions~\ref{def:kcorr} and~\ref{def:n-corr} to describe the actions of $N_n$, the Galois group $\Gal(\Q(\omega)\mid\Q)$, and $\Lin(P_n)$ on $\Irr(P_n)$.

\begin{notation}
As in Notation~\ref{not:lincp}, we write $\Irr(P_p)=\{\phi_{\varepsilon}\mid\varepsilon\in [0,p-1]\}$. We extend the action of $\sigma_1^{(1)}$ and its inverse $\tau$ to the set $[0,p]$ by fixing both $0$ and $p$. Thus, 
\[
(0)\sigma_1^{(1)}=(0)\tau=0, \qquad (p)\sigma_1^{(1)}=(p)\tau=p.
\]
\end{notation}

\subsection{Normalizer action on irreducible characters of $P_n$}\label{sec:nor} 
The action of $N_{n}$ on $\Irr(P_n)$ is given by 
\begin{equation*}
\chi^{h}(g)=\chi(g^{h^{-1}})=\chi(hgh^{-1}),
\qquad \chi\in\Irr(P_n), g\in P_n, h\in N_n.
\end{equation*}
It is convenient to study instead the action of $N_n$ on $\cF_n$, defined by
\begin{equation*}
\cT^{h}:=\Phi_n(\Psi_n(\cT)^h), \qquad \cT\in\cF_n, h\in N_n.
\end{equation*}
These actions are isomorphic. In particular, for all $\phi\in\Irr(P_n)$, $\mathcal{S}\in\cF_n$, and $h\in N_n$, we can write
\begin{equation}\label{eq:azequivnor}
\begin{split}
&\cT(\phi)^h=\cT(\phi^h), \\ 
&\theta(\mathcal{S}^h)=\theta(\mathcal{S})^h. 
\end{split}
\end{equation}
Since characters are class functions, $P_n$ acts trivially, and it suffices to study the induced action of $N_n/ P_n$. 

We start with $n=p^k$ for some integer $k\geq 1$. From Definition~\ref{def:N_p^k}, $\langle \sigma_j^{(k)}P_{p^k} \mid j\in [1,k] \rangle=N_{p^k}/P_{p^k}$. Therefore, we represent $\cT$-functions as in Notation~\ref{not:T=()}, and we focus on the generators $\{\sigma_j^{(k)}\mid j\in [1,k]\}$.

\begin{theorem}\label{th:1}
Let $p$ be a prime and $k\geq 1$ an integer. Let
$\cT\in\cF_{p^k}$ be written as $\cT=(\cT_1 \mid \dots \mid \cT_p; x)$. Then, 
\[ \cT^{\sigma_j^{(k)}}=
\begin{cases}
(\cT_1^{\sigma_j^{(k-1)}} \mid \dots \mid \cT_p^{\sigma_j^{(k-1)}}; x) & \text{if } j<k \\
(\cT_{(1)\tau} \mid \dots \mid \cT_{(p)\tau}; (x)\tau) & \text{if } j=k.
\end{cases}
\]
\end{theorem}

\begin{proof}
We proceed by induction on $k$. If $k=1$, then $j=1$ and $\sigma_j^{(k)}=\sigma_1^{(1)}$. By Notation~\ref{not:cTfunc} and Definition~\ref{def:kcorr}, we have $\cT=(x)=\cT^{x}=\Phi_p(\phi_x)=\cT(\phi_x)$ for some $x\in[0,p-1]$. Consequently,
\[
\cT^{\sigma_1^{(1)}}=\cT(\phi_{x})^{\sigma_1^{(1)}}=\cT((\phi_{x})^{\sigma_1^{(1)}})=\cT(\phi_{(x)\tau})=\cT^{(x)\tau}=((x)\tau),
\]
where the second equality follows from \eqref{eq:azequivnor}, and the third from \eqref{eq:t}. This establishes the base case. 
We now suppose that $k>1$ and distinguish cases according to the values of $x$ and $j$. In each case, we denote by $\cR$ the $\cT$-function prescribed by the statement, and we prove that $\cT^{\sigma_j^{(k)}}=\cR$ by showing that $\theta(\cT^{\sigma_j^{(k)}})=\theta(\cR)$. This suffices, by Lemma~\ref{lem:Tbij}. In what follows, we repeatedly use Lemmas~\ref{lem:clif} and~\ref{lem:Tbij}, Proposition~\ref{prop:sigmasuP} and the formula for values of characters of wreath products given in \cite[4.4.10]{JK}. To avoid redundancy, we cite these results only once. We let $g=(f;\delta)=(f_1,\dotsc,f_p;\delta)\in P_{p^{k-1}}\wr P_p=P_{p^k}$, where $f=(f_1,\dotsc,f_p)\in B=(P_{p^{k-1}})^{\times p}$ and $\delta\in P_p$. We verify that $\theta(\cT^{\sigma_j^{(k)}})(g)=\theta(\cR)(g)$ in each case.

\noindent\textbf{\textit{\underline{Case $x\neq p$.}}}  Since $\cT\in\cF_{p^k}$, it follows that $\cT_1=\cdots=\cT_p=\Phi_{p^{k-1}}(\varphi)=\cT(\varphi)$ for some $\varphi\in\Irr(P_{p^{k-1}})$. Then, $\cT=(\cT(\varphi) \mid \dots \mid \cT(\varphi);x)$ and $\theta(\cT)=\cX(\varphi;x)$, by Lemma~\ref{lem:Tbij}. Thus, using $\theta(\cT^{\sigma_j^{(k)}})=\theta(\cT)^{\sigma_j^{(k)}}$, 
\begin{equation*}
\theta(\cT^{\sigma_j^{(k)}})=\cX(\varphi;x)^{\sigma_j^{(k)}}.
\end{equation*}

\noindent\textbf{Subcase $j<k$.}	Let $\cR=(\cT_1^{\sigma_j^{(k-1)}} \mid \dots \mid \cT_1^{\sigma_j^{(k-1)}}; x)$, with $\cT_1=\cT(\varphi)$. Since $\cT(\varphi)^{\sigma_j^{(k-1)}}=\cT(\varphi^{\sigma_j^{(k-1)}})$, then $\cR=(\cT(\varphi^{\sigma_j^{(k-1)}})\mid\dots\mid\cT(\varphi^{\sigma_j^{(k-1)}});x)$, and $\theta(\cR)=\cX(\varphi^{\sigma_j^{(k-1)}};x)$. By Proposition~\ref{prop:sigmasuP}, we have that
\begin{equation}\label{eq:ast}
\theta(\cT^{\sigma_j^{(k)}})(g)=\cX(\varphi;x)^{\sigma_j^{(k)}}(g)=\cX(\varphi;x)((f;\delta)^{(\sigma_j^{(k)})^{-1}})=\cX(\varphi;x)((f_1^{(\sigma_j^{(k-1)})^{-1}},\dotsc,f_p^{(\sigma_j^{(k-1)})^{-1}};\delta)).
\end{equation}
If $\delta=1$, then $g,g^{(\sigma_j^{(k)})^{-1}}\in B$. By Lemma~\ref{lem:clif}, $\cX(\varphi;x)\down_B=\varphi^{\times p}$ and $\cX(\varphi^{\sigma_j^{(k-1)}};x)\down_B=(\varphi^{\sigma_j^{(k-1)}})^{\times p}$, so
\begin{align*}
\theta(\cT^{\sigma_j^{(k)}})(g)&=\cX(\varphi;x)((f_1^{(\sigma_j^{(k-1)})^{-1}},\dotsc,f_p^{(\sigma_j^{(k-1)})^{-1}};1))=\prod_{i=1}^{p}\varphi(f_i^{(\sigma_j^{(k-1)})^{-1}})=\prod_{i=1}^{p}\varphi^{\sigma_j^{(k-1)}}(f_i) \\
&=\cX(\varphi^{\sigma_j^{(k-1)}};x)((f_1,\dotsc,f_p;1))=\theta(\cR)(g).
\end{align*}
If $\delta\neq 1$ instead, we use Lemma \cite[4.4.10]{JK} in \eqref{eq:ast} to obtain
\begin{align*}
\theta(\cT^{\sigma_j^{(k)}})(g)&=\phi_{x}(\delta)\cdot\varphi(f_1^{(\sigma_j^{(k-1)})^{-1}}f_{(1)\delta^{-1}}^{(\sigma_j^{(k-1)})^{-1}} \cdots f_{(1)\delta^{-p+1}}^{(\sigma_j^{(k-1)})^{-1}})=\phi_{x}(\delta)\cdot\varphi((f_1f_{(1)\delta^{-1}} \cdots f_{(1)\delta^{-p+1}})^{(\sigma_j^{(k-1)})^{-1}}) \\
&=\phi_{x}(\delta)\cdot\varphi^{\sigma_j^{(k-1)}}(f_1f_{(1)\delta^{-1}} \cdots f_{(1)\delta^{-p+1}})=\cX(\varphi^{\sigma_j^{(k-1)}};x)((f;\delta))=\theta(\cR)(g).
\end{align*}

\noindent\textbf{Subcase $j=k$.}	Let $\cR=(\cT_1\mid\dots\mid\cT_1;(x)\tau)$, where $\cT_1=\cT(\varphi)$ and $x\neq p$. As a consequence, $(x)\tau\neq p$. Hence,
$\cR=(\cT(\varphi)\mid\dots\mid\cT(\varphi);(x)\tau)$ and $\theta(\cR)=\cX(\varphi;(x)\tau)$. We evaluate
\begin{align}\label{eq:astast}
\theta(\cT^{\sigma_k^{(k)}})(g)&=\cX(\varphi;x)^{\sigma_k^{(k)}}(g)=\cX(\varphi;x)((f;\delta)^{(\sigma_k^{(k)})^{-1}})=\cX(\varphi;x)((f_{(1)\sigma_1^{(1)}},\dotsc,f_{(p)\sigma_1^{(1)}};\delta^{\tau})).
\end{align}
If $\delta=1$, then $g,g^{(\sigma_k^{(k)})^{-1}}\in B$. From $\cX(\varphi;x)\down_B=\varphi^{\times p}=\cX(\varphi;(x)\tau)\down_B$, we derive that
\begin{align*}
\theta(\cT^{\sigma_k^{(k)}})(g)&=\cX(\varphi;x)((f_{(1)\sigma_1^{(1)}},\dotsc,f_{(p)\sigma_1^{(1)}};1)=\prod_{i=1}^{p}\varphi(f_{(i)\sigma_1^{(1)}})=\prod_{i=1}^{p}\varphi(f_i) \\
&=\cX(\varphi;(x)\tau)((f_1,\dotsc,f_p;1))=\theta(\cR)(g).
\end{align*}
On the other hand, let $\delta\neq 1$. From \eqref{eq:astast}, we compute
\begin{align*}
\theta(\cT^{\sigma_k^{(k)}})(g)&=\phi_{x}(\delta^{\tau})\cdot\varphi(f_{(1)\sigma_1^{(1)}}f_{((1)\sigma_1^{(1)}) \delta^{-1}} \cdots f_{((1)\sigma_1^{(1)}) \delta^{-p+1}})=(\phi_{x})^{\sigma_1^{(1)}}(\delta)\cdot\varphi(\prod_{i=1}^{p}f_{((1)\sigma_1^{(1)}) \delta^{-(i-1)}})\\
&=\phi_{(x)\tau}(\delta)\cdot\varphi(\prod_{i=1}^{p}f_{((1)\sigma_1^{(1)}) \delta^{-(i-1)}})=\phi_{(x)\tau}(\delta)\cdot\varphi(\prod_{i=1}^{p}f_{(1)\delta^{-(i-1)}}) \\
&=\phi_{(x)\tau}(\delta)\cdot\varphi(f_1f_{(1)\delta^{-1}} \cdots f_{(1)\delta^{-p+1}})=\cX(\varphi;(x)\tau)((f;\delta))=\theta(\cR)(g).
\end{align*} 
For the third equality we use \eqref{eq:t}, while  the fourth follows from \cite[4.2.5]{JK}, which states that the two arguments of $\varphi$ lie in the same $P_{p^{k-1}}$-conjugacy class. This completes the case $x\neq p$. 

\noindent\textbf{\textit{\underline{Case $x=p$.}}}  We have that $\cT_1=\Phi_{p^{k-1}}(\varphi_1),\dotsc,\cT_p=\Phi_{p^{k-1}}(\varphi_p)$ for some mixed $\varphi_1,\dotsc,\varphi_p\in\Irr(P_{p^{k-1}})$. Thus, $\cT=(\cT(\varphi_1)\mid\dots\mid\cT(\varphi_p);p)$ and $\theta(\cT)=\varphi_1\times\cdots\times\varphi_p\up^{P_{p^k}}$, with $\varphi_1\times\cdots\times\varphi_p\in\Irr(B)$ and $B\trianglelefteq P_{p^k}$. As shown in the proof of \cite[Lemma 4.2]{OlssonPaper}, $B= \langle y\in P_{p^k} \mid y \text{ has a fixed point}  \rangle\trianglelefteq N_{p^k}$. Therefore, $N_{p^k}$ acts by conjugation on $\Irr(B)$. Since $\theta(\cT)\down_B=\sum_{\pi\in P_p}\varphi_{(1)\pi}\times\cdots\times\varphi_{(p)\pi}$, we obtain
\begin{align*}
1&=[\theta(\cT)\down_B,\varphi_1\times\cdots\times\varphi_p]=[\theta(\cT)\down_B{}^{\sigma_j^{(k)}},(\varphi_1\times\cdots\varphi_p)^{\sigma_j^{(k)}}]
=[\theta(\cT)^{\sigma_j^{(k)}}\down_B,(\varphi_1\times\cdots\times\varphi_p)
^{\sigma_j^{(k)}}] \\
&=[\theta(\cT^{\sigma_j^{(k)}})\down_B,(\varphi_1\times\cdots\times\varphi_p)
^{\sigma_j^{(k)}}]=[\theta(\cT^{\sigma_j^{(k)}}),(\varphi_1\times\cdots\varphi_p)^{\sigma_j^{(k)}}\up^{P_{p^{k}}}].
\end{align*}
Hence, $\theta(\cT^{\sigma_j^{(k)}})$ is an irreducible constituent of $(\varphi_1\times\cdots\times\varphi_p)^{\sigma_j^{(k)}}\up^{P_{p^k}}$, occurring with multiplicity one. As the degrees of these characters coincide, we conclude that 
\[
\theta(\cT^{\sigma_j^{(k)}})=(\varphi_1\times\cdots\times\varphi_p)^{\sigma_j^{(k)}}\up^{P_{p^k}}.
\] 

\noindent\textbf{Subcase $j<k$.} Let $\cR=(\cT_1^{\sigma_j^{(k-1)}}\mid\dots\mid\cT_p^{\sigma_j^{(k-1)}};p)$, with $\cT_i^{\sigma_j^{(k-1)}}=\cT(\varphi_i)^{\sigma_j^{(k-1)}}=\cT(\varphi_i^{\sigma_j^{(k-1)}})$ for all $i$. This implies that $\cR=(\cT(\varphi_1^{\sigma_j^{(k-1)}})\mid\dots\mid\cT(\varphi_p^{\sigma_j^{(k-1)}});p)$  and $\theta(\cR)=\varphi_1^{\sigma_j^{(k-1)}}\times\cdots\times\varphi_p^{\sigma_j^{(k-1)}}\up^{P_{p^k}}$. If $\delta=1$,
\begin{align*}
&(\varphi_1\times\cdots\times\varphi_p)^{\sigma_j^{(k)}}(g)=
\varphi_1\times\cdots\times\varphi_p((f;1)^{(\sigma_j^{(k)})^{-1}})=
\varphi_1\times\cdots\times\varphi_p((f_1^{(\sigma_j^{(k-1)})^{-1}},\dotsc,f_p^{(\sigma_j^{(k-1)})^{-1}};1)) \\
&=\prod_{i=1}^{p}
\varphi_i(f_i^{(\sigma_j^{(k-1)})^{-1}})
=\prod_{i=1}^{p}\varphi_i^{\sigma_j^{(k-1)}}(f_i)=\varphi_1^{\sigma_j^{(k-1)}}\times\cdots\times\varphi_p^{\sigma_j^{(k-1)}}((f;1))=\varphi_1^{\sigma_j^{(k-1)}}\times\cdots\times\varphi_p^{\sigma_j^{(k-1)}}(g).
\end{align*}
Thus $(\varphi_1\times\cdots\times\varphi_p)^{\sigma_j^{(k)}}=\varphi_1^{\sigma_j^{(k-1)}}\times\cdots\times\varphi_p^{\sigma_j^{(k-1)}}$, and inducing from both sides we get $\theta(\cT^{\sigma_j^{(k)}})=\theta(\cR)$.

\noindent\textbf{Subcase $j=k$.}	Let $\cR=(\cT_{(1)\tau}\mid\dots\mid\cT_{(p)\tau}; p)$ and $\theta(\cR)=\varphi_{(1)\tau}\times\cdots\times\varphi_{(p)\tau}\up^{P_{p^k}}$. If $\delta=1$, 
\begin{align*}
&(\varphi_1\times\cdots\times\varphi_p)^{\sigma_k^{(k)}}(g)=\varphi_1\times\cdots\times\varphi_p((f;1)^{(\sigma_k^{(k)})^{-1}})=
\varphi_1\times\cdots\times\varphi_p((f_{(1)\sigma_1^{(1)}},\dotsc,f_{(p)\sigma_1^{(1)}};1)) \\
&=\prod_{i=1}^{p}\varphi_i(f_{(i)\sigma_1^{(1)}})=\prod_{j=1}^{p}\varphi_{(j)\tau}(f_j)=\varphi_{(1)\tau}\times\cdots\times\varphi_{(p)\tau}((f;1))=\varphi_{(1)\tau}\times\cdots\times\varphi_{(p)\tau}(g).
\end{align*}
Hence $(\varphi_1\times\cdots\times\varphi_p)^{\sigma_k^{(k)}}=\theta_{(1)\tau}\times\cdots\times\theta_{(p)\tau}$, and inducing from both sides we obtain $\theta(\cT^{\sigma_k^{(k)}})=\theta(\cR)$. This completes the proof.
\end{proof}

From this result, we derive a straightforward method to construct a $p^k$-labeling function representing $\cT^{\sigma_j^{(k)}}$, provided a representative of $\cT$ is known.

\begin{corollary}\label{cor:1}
Let $p$ be a prime and $k\geq 1$ an integer. Let $\cT\in\cF_{p^k}$ and $\ts=s_1\dots s_{\ell(\ts)} \in\sk_{p^k}$. Then,

\[
\cT^{\sigma_j^{(k)}}(\ts) =
\left\{
\begin{array}{ll}
  \begin{cases}
    (\cT(\emptyset))\tau & \text{if } \ell(\ts) = 0, \\
    \cT((s_1)\tau \mid {}^{1}\ts) & \text{if } \ell(\ts) \neq 0
  \end{cases} & \text{if } j = k, \\[1.5em]
  
  \begin{cases}
    \cT(\ts) & \text{if } \ell(\ts) < k - j, \\
    (\cT(\ts))\tau & \text{if } \ell(\ts) = k - j, \\
    \cT(\ts^{k-j} \mid (s_{k-j+1})\tau \mid {}^{k-j+1}\ts) & \text{if } \ell(\ts) > k - j
  \end{cases} & \text{if } j < k.
\end{array}
\right.
\]
\end{corollary}

\begin{proof}
The result is an immediate consequence of Theorem~\ref{th:1}, since $(\cT((s_1)\tau))^{\giu}(^{1}\ts) = \cT((s_1)\tau \mid {}^{1}\ts)$ 
and $(\cT(\ts^{k-j} \mid (s_{k-j+1})\tau))^{\giu}({}^{k-j+1}\ts) = \cT(\ts^{k-j} \mid (s_{k-j+1})\tau \mid {}^{k-j+1}\ts)$. 
\end{proof}

Let $n\in\N$. We describe the action of the generators of $N_n$ (see Definition~\ref{def:N_n}). Since $\cT^{g(i,j,\ell)}=\cT$ for every $\cT\in\cF_n$ and $(i,j,\ell)\in L(n)$, by Corollary~\ref{cor:1} and Proposition~\ref{prop:NsuPgenerale} we obtain the next theorem. 

\begin{theorem}\label{th:2}
Let $p$ be a prime and $n\in\N$ with $p$-adic expansion as in Notation~\ref{not:p-adicn}.
Let $\cT\in\cF_n$, $\sigma(i,j,\ell)\in N_n$ and $\rho_y\in \langle \xi_y, \zeta_y \rangle$, for $(i,j,\ell)\in L(n)$ and $y\in[1,t]$. Then, for every $\mathbf{\ts}\in\sk_n$,

\[
\cT^{\sigma(i,j,\ell)}(\mathbf{\ts}) = \left( \cT_1(\mathbf{\ts}_1),\dotsc, \cT_{q_{(i-1)}+1}(\mathbf{\ts}_{q_{(i-1)}+1}), \dotsc, \cT_{q_{(i-1)}+j}^{\sigma_\ell^{(k_i)}}(\mathbf{\ts}_{q_{(i-1)}+j}), \dotsc, \cT_{q_t}(\mathbf{\ts}_{q_t}) \right),
\]

\[
\cT^{\rho_y}(\mathbf{\ts})= \left( \cT_1(\mathbf{\ts}_1),\dotsc, \cT_{q_{(y-1)}+(1)\rho_y^{-1}}(\mathbf{\ts}_{q_{(y-1)}+(1)\rho_y^{-1}}),\dotsc, \cT_{q_{(y-1)}+(a_y)\rho_y^{-1}}(\mathbf{\ts}_{q_{(y-1)}+(a_y)\rho_y^{-1}}),\dotsc, \cT_{q_t}(\mathbf{\ts}_{q_t}) \right)
.
\]
\end{theorem}

\begin{remark}
These results extend the study of the $N_n$-action on $\Lin(P_n)$ in \cite{law,gll,GANT} to the full set of irreducible characters $\Irr(P_n)$. Thus, the linear case is recovered by restricting to $\Lin(P_n)$. More precisely, let $n\in\N$ have $p$-adic expansion as in Notation~\ref{not:p-adicn}. Any $\lambda\in\Lin(P_n)$ decomposes as a direct product of linear characters of the factors $P_{p^{\mathtt{o}(j)}}$, for $j\in [1,q_t]$. For $\mu\in\Lin(P_{p^k})$ and $\ts\in\sk_{p^k}$, we have that $\mathcal{L}=\cT(\mu)$ satisfies $
\mathcal{L}(\ts)=\mathcal{L}((s_1)\gamma_1(s_2)\gamma_2 \dots (s_{\ell(\ts)})\gamma_{\ell(\ts)})$
for every $\gamma_i\in P_p$ with $i\in [1,\ell(\ts)]$. Equivalently, $\mathcal{L}$ takes values in $[0,p-1]$ and depends only on $\ell(\ts)$, being constant on sequences of equal length. Hence $\mathcal{L}$ has a unique representative, as in Example~\ref{ex:fromttoT} and Figure~\ref{fig:linear}. Let $\ts_i\in\sk_{p^k}$ with $\ell(\ts_i)=i-1$ for $i\in [1,k]$. The map $\mu\mapsto\mathcal{L}\mapsto(\mathcal{L}(\ts_1),\dotsc,\mathcal{L}(\ts_k))$ defines a bijection between $\Lin(P_{p^k})$ and $[0,p-1]^{[1,k]}$. This is compatible with the parametrization of linear characters in \cite[Section 2]{GANT}.
\end{remark}

\subsection{Galois action on irreducible characters of $P_n$}\label{sec:gal}
Let $\omega=e^{2\pi i/p}$ and write $\cG=\Gal(\Q(\omega)\mid\Q)$. By \cite[Theorem 4.4.8]{JK}, $\Q(\omega)$ is a splitting field for $P_n$. The group $\cG$ acts naturally on $\Irr(P_n)$ via
\[
\chi^{\rho}(g)=(\chi(g))^{\rho}, \qquad \chi\in\Irr(P_n), g\in P_n, \rho\in\cG.
\]
We study the action of $\cG$ on $\cF_n$ defined by
\[
\cT^{\rho}=\Phi_n(\Psi_n(\cT)^{\rho}), \qquad \cT\in\cF_n, \rho\in\cG,
\]
as these two actions are isomorphic. In particular, for all $\phi\in\Irr(P_n),\mathcal{S}\in\cF_n$, and $\rho\in\cG$, we have
\begin{equation}\label{eq:azgaleq}
\begin{split}
&\cT(\phi)^{\rho}=\cT(\phi^{\rho}), \\
&\theta(\mathcal{S}^{\rho})=\theta(\mathcal{S})^{\rho}.
\end{split}
\end{equation}

\begin{definition}\label{def:sigma}
Let $\sigma$ be the generator of $\cG$ satisfying $(\phi_{\varepsilon})^{\sigma}=(\phi_{\varepsilon})^{\sigma_1^{(1)}}=\phi_{(\varepsilon)\tau}$ for all $\varepsilon\in [0,p-1]$. 
\end{definition}

As shown in Notation~\ref{not:lincp}, $\sigma$ is well defined. Since $\cG=\langle \sigma \rangle$, it suffices to determine the action of $\sigma$ on $\cF_n$ in order to describe the action of $\cG$. A combinatorial description is provided in Theorems~\ref{th:3} and~\ref{th:4}. Specifically, we construct a labeling function associated with $\cT^{\sigma}$ from a given representative of $\cT$.

\begin{theorem}\label{th:3}
Let $p$ be a prime and $k\geq 1$ an integer. Let $\cT\in\cF_{p^k}$ and $\ts\in\sk_{p^k}$. Then,
$\cT^{\sigma}(\ts)=(\cT(\ts))\tau$.
\end{theorem}

\begin{proof}
We proceed by induction on $k$. If $k=1$, then $\ts=\emptyset$. From Definition~\ref{def:sigma}, Lemma~\ref{lem:Tbij} and equations \eqref{eq:azequivnor}, \eqref{eq:azgaleq}, it follows that $\cT^{\sigma}=\cT^{\sigma_1^{(1)}}$. Thus, the base case holds by applying Corollary~\ref{cor:1}. Let $k\ge 1$. Let $\cT=(\cT_1\mid\dots\mid\cT_p;x)$ as in Notation~\ref{not:T=()}. We consider the two following different cases.

\noindent\underline{\textbf{\textit{Case \textnormal{$x\neq p$}.}}} 
From Definition~\ref{def:admT} and Lemma~\ref{lem:Tbij}, we deduce that there exists some $\varphi\in\Irr(P_{p^{k-1}})$ such that $\cT_1=\dots=\cT_p=\cT(\varphi)\in\cF_{p^{k-1}}$. Thus, $\cT=(\cT(\varphi)\mid\dots\mid\cT(\varphi);x)$ and $\theta(\cT)=\cX(\varphi;x)$. Similarly, since $(x)\tau\neq p$, we let $\cR=(\cT(\varphi^{\sigma})\mid\dots\mid \cT(\varphi^{\sigma}); (x)\tau)\in\cF_{p^k}$ and $\theta(\cR)=\cX(\varphi^{\sigma};(x)\tau)$. Let $B=(P_{p^{k-1}})^{\times p}$ denote the base group of $P_{p^{k-1}}\wr P_p=P_{p^k}$. Using \eqref{eq:azgaleq} and Lemma~\ref{lem:clif}, we deduce that
\[
\theta(\cT^{\sigma})\down_B=\theta(\cT)^{\sigma}\down_B=(\theta(\cT)\down_B)^{\sigma}=(\cX(\varphi;x)\down_B)^{\sigma}=(\varphi^{\times p})^{\sigma}=(\varphi^{\sigma})^{\times p}=\cX(\varphi^{\sigma};(x)\tau)\down_B=\theta(\cR)\down_B.
\]
Now let $g=(f;\delta)\in P_{p^k}\smallsetminus B$, and set $\rho(f;\delta)=\prod_{i=1}^{p}f_{(1)\delta^{-i+1}}$. Using \cite[4.4.10]{JK}, we compute
\begin{align*}
\theta(\cT^{\sigma})(g)&=\theta(\cT)^{\sigma}(g)=
(\cX(\varphi;x)(f;\delta))^{\sigma}=(\phi_{x}(\delta)\cdot\varphi(\rho(f;\delta)))^{\sigma}=(\phi_{x}(\delta))^{\sigma}\cdot(\varphi(\rho(f;\delta)))^{\sigma} \\
&=\phi_{(x)\tau}(\delta)\cdot\varphi^{\sigma}(\rho(f;\delta))=\cX(\varphi^{\sigma};(x)\tau)(f;\delta)=\theta(\cR)(g).
\end{align*}
This implies that $\theta(\cT^{\sigma})=\theta(\cR)$, and thus $\cT^{\sigma}=\cR$. Note that
$\cT^{\sigma}(\emptyset)=\cR(\emptyset)=(x)\tau=(\cT(\emptyset))\tau$. If $\ts\neq\emptyset$, 
\begin{align*}
\cT^{\sigma}(\ts)&=\cR(\ts)=\cR(s_1| ^{1}\ts)=(\cR({s_1}))^{\giu}(^{1}\ts)=\cT(\varphi^{\sigma})(^{1}\ts)=\cT(\varphi)^{\sigma}(^{1}\ts)=(\cT(\varphi)(^{1}\ts))\tau=((\cT(s_1))^{\giu}(^{1}\ts))\tau \\
&=(\cT(s_1| ^{1}\ts))\tau=(\cT(\ts))\tau,
\end{align*}
where the sixth equality holds by the inductive hypothesis. The theorem follows in this instance. 

\noindent\underline{\textbf{\textit{Case \textnormal{$x=p$}.}}} 
We have that $\cT_1=\Phi_{p^k}(\varphi_1),\dotsc,\cT_p=\Phi_{p^k}(\varphi_p)$ for some mixed $\varphi_1,\dotsc,\varphi_p\in\Irr(P_{p^{k-1}})$. Thus, by Lemma~\ref{lem:Tbij}, $\cT=(\cT(\varphi_1)\mid\dots\mid\cT(\varphi_p);p)$ and $\theta(\cT)=\varphi_1\times\cdots\times\varphi_p
\up^{P_{p^k}}$. From \eqref{eq:azgaleq} and the definition of induced character, it follows that $\theta(\cT^{\sigma})=\theta(\cT)^{\sigma}=(\varphi_1\times\cdots\times\varphi_p
\up^{P_{p^k}})^{\sigma}=\varphi_1^{\sigma}\times\cdots\times\varphi_p^{\sigma}\up^{P_{p^k}}$. Hence, we derive that $\cT^{\sigma}=(\cT(\varphi_1^{\sigma})\mid\dots\mid\cT(\varphi_p^{\sigma});p)$. Observe that $\cT^{\sigma}(\emptyset)=p=(\cT(\emptyset))\tau$. If $\ts\neq \emptyset$, then 
\begin{align*}
\cT^{\sigma}(\ts)&=\cT^{\sigma}(s_1|^{1}\ts)
=(\cT^{\sigma}(s_1))^{\giu}(^{1}\ts) =\cT(\varphi_{s_1}^{\sigma})(^{1}\ts)=\cT(\varphi_{s_1})^{\sigma
}(^{1}\ts)=(\cT(\varphi_{s_1})(^{1}\ts))\tau=((\cT(s_1))^{\giu}(^{1}\ts))\tau \\
&=(\cT(s_1|^{1}\ts))\tau=(\cT(\ts))\tau,
\end{align*}
where the fifth equality holds by the inductive hypothesis. This completes the proof.
\end{proof}

\begin{theorem}\label{th:4}
Let $p$ be a prime and $n\in\N$ with $p$-adic expansion as in Notation~\ref{not:p-adicn}. Let $\cT\in\cF_n$. Then, $\cT^{\sigma}(\mathbf{\ts})=\left( (\cT_1(\mathbf{\ts}_1))\tau,\dotsc, (\cT_{q_t}(\mathbf{\ts}_{q_t}))\tau \right)$ for every $\mathbf{\ts}\in\sk_n$.
\end{theorem}

\begin{proof}
By Definitions~\ref{def:cTfuncgen} and~\ref{def:n-corr}, if $\cT=(\cT_1,\dotsc,\cT_{q_t})$ then $\theta(\cT)=\theta(\cT_1)\times\cdots\times\theta(\cT_{q_t})$. From
\[
\theta(\cT^{\sigma})=\theta(\cT)^{\sigma}=(\theta(\cT_1)\times\cdots\times\theta(\cT_{q_t}))^{\sigma}=\theta(\cT_1)^{\sigma}\times\cdots\times\theta(\cT_{q_t})^{\sigma}=\theta(\cT_1^{\sigma})\times\cdots\times\theta(\cT_{q_t}^{\sigma}),
\]
we deduce that $\cT^{\sigma}=(\cT_1^{\sigma},\dotsc,\cT_{q_t}^{\sigma})$. The statement now follows immediately from Theorem~\ref{th:3}.
\end{proof}

Let $p>2$. The unique nontrivial involution of $\cG=\langle\sigma\rangle\cong C_{p-1}$ is $\sigma^{\tfrac{p-1}{2}}$, which is the field automorphism $\overline{\cdot}: \Q(\omega) \to \Q(\omega)$ given by complex conjugation. From Theorem~\ref{th:3} and 
\[
\tau^{\tfrac{p-1}{2}}=(1 \quad p-1)(2 \quad p-2)\cdots(\tfrac{p-1}{2} \quad \tfrac{p+1}{2})=\prod_{i=1}^{(p-1)/2}(i \quad p-i),
\] 
we obtain the following. The general case is obtained by the same argument as in Theorem~\ref{th:4}, using Corollary~\ref{cor:6}. 

\begin{corollary}\label{cor:6}
Let $p$ be a prime and $k\geq 1$ an integer. Let $\cT\in\cF_{p^k}$ and $\ts\in\sk_{p^k}$. Then,  
\[
\overline{\cT}(\ts)=(\cT(\ts))\tau^{\tfrac{p-1}{2}}=
\begin{cases}
p - \cT(\ts) & \ \text{if } p\nmid\cT(\ts) \\
\cT(\ts) & \ \text{if } p\mid\cT(\ts)
\end{cases}
\]
\end{corollary}

\begin{remark}
A character is \textit{rational} if it takes values in $\Q$. If $\phi\in\Irr(P_n)$ is rational, then $\overline{\phi}=\phi$. When $p>2$, the only rational irreducible character of $P_n$ is the trivial one, since by Definitions~\ref{def:admT} and~\ref{def:cTfuncgen} and Corollary~\ref{cor:6}, it follows that $\overline{\phi}\neq\phi$ for every $\phi\in\Irr(P_n)\smallsetminus\{\triv_{P_n}\}$. If $p=2$, we have $\Q(\omega)=\Q$ and hence every irreducible character of $P_n$ is rational. 
\end{remark}

\subsection{Multiplication by linear characters of $P_n$}
Let $n\in\N$ have $p$-adic expansion as in Notation~\ref{not:p-adicn}. Consider $\varphi=\varphi_1\times\cdots\times\varphi_{q_t}\in\Irr(P_n)$ and $\lambda=\lambda_1\times\cdots\times\lambda_{q_t}\in\Lin(P_n)$. The product $\varphi\cdot\lambda\in\Irr(P_n)$ is defined by $(\varphi\cdot\lambda)(g)=\varphi(g)\lambda(g)$ for all $g\in P_n$, and $\Lin(P_n)$ acts by multiplication on $\Irr(P_n)$.  Since $\varphi\cdot\lambda=(\varphi_1\cdot\lambda_1)\times\cdots\times(\varphi_{q_t}\cdot\lambda_{q_t})$, it suffices to consider the case $n=p^k$ for some integer $k\geq 1$. Through the $\cT$-correspondence, this action admits a simple combinatorial description.

\begin{theorem}\label{th:5}
Let $p$ be a prime and $k\geq 1$ an integer. For every $a\in\N$, let $\underline{a}$ be the unique element in $[0,p-1]$ such that $a \equiv \underline{a} \bmod p$. Let $\varphi\in\Irr(P_{p^k})$, $\lambda\in\Lin(P_{p^k})$ and $\ts\in\sk_{p^k}$. Then,
\[
\cT(\varphi\cdot\lambda)(\ts)=
\begin{cases}
\cT(\varphi)(\ts) & \ \text{if } \cT(\varphi)(\ts)=p, \\
\underline{\cT(\varphi)(\ts) + \cT(\lambda)(\ts)} & \ \text{if } \cT(\varphi)(\ts)\neq p.
\end{cases}
\]
\end{theorem}

\begin{proof}
We proceed by induction on $k$. If $k=1$, then $\ts=\emptyset$ and there exist $\varepsilon_1,\varepsilon_2\in [0,p-1]$ such that $\varphi=\phi_{\varepsilon_1}$ and $\lambda=\phi_{\varepsilon_2}$. Let $g\in P_p$ and compute
\begin{align*}
(\phi_{\varepsilon_1}\cdot\phi_{\varepsilon_2})(g)=\phi_{\varepsilon_1}(g)\phi_{\varepsilon_2}(g)=(\phi_1(g))^{\varepsilon_1}(\phi_1(g))^{\varepsilon_2}=(\phi_1(g))^{\underline{\varepsilon_1 + \varepsilon_2}}=\phi_{\underline{\varepsilon_1+\varepsilon_2}}(g).
\end{align*}
This implies that $\cT(\phi_{\varepsilon_1}\cdot\phi_{\varepsilon_2})=\cT(\phi_{\underline{\varepsilon_1+\varepsilon_2}})$, and therefore $\cT(\phi_{\varepsilon_1}\cdot\phi_{\varepsilon_2})(\emptyset)=\underline{\cT(\phi_{\varepsilon_1})(\emptyset)+\cT(\phi_{\varepsilon_2})(\emptyset)}$ holds. Suppose $k>1$. Let $g=(f;\delta)=(f_1,\dotsc,f_p;\delta)\in P_{p^{k-1}}\wr P_p=P_{p^k}$, and let $B=(P_{p^{k-1}})^{\times p}$. By Lemma~\ref{lem:chP-trees}, 
\[
\lambda=\cX(\mu;y)
\]
for some $\mu\in\Lin(P_{p^{k-1}})$ and $y\in [0,p-1]$, since $\lambda\in\Lin(P_{p^k})$. Using Notation~\ref{not:T=()} and Lemma~\ref{lem:Tbij}, 
\[
\cT(\lambda)=(\cT(\mu)\mid\dots\mid\cT(\mu);y), \qquad \cT(\varphi)=(\cT_1\mid\dots\mid\cT_p;x)
\] 
for some $\cT_1,\dotsc,\cT_p\in\cF_{p^{k-1}}$ and $x\in [0,p]$. We consider cases separately, according to the value of $x$. 

\noindent\underline{\textbf{\textit{Case \textnormal{$x\neq p$}.}}} From Definition~\ref{def:admT} and Lemma~\ref{lem:Tbij}, we deduce that there exists some $\psi\in\Irr(P_{p^{k-1}})$ such that $\cT_1=\ldots=\cT_p=\cT(\psi)\in\cF_{p^{k-1}}$. As a consequence,
\[
\cT(\varphi)=(\cT(\psi)\mid\dots\mid\cT(\psi);x), \qquad \varphi=\cX(\psi;x).
\]
Let $\mathcal{R}=(\cT(\psi\cdot\mu) \mid \dots \mid \cT(\psi\cdot\mu); \underline{x+y})\in\cF_{p^k}$,
and $\theta(\cR)=\cX(\psi\cdot\mu; \underline{x+y})$. By Lemma~\ref{lem:clif}, we derive
\[
(\varphi\cdot\lambda)\down_B=\varphi\down_B\cdot\lambda\down_B=\cX(\psi;x)\down_B\cdot\cX(\mu;y)\down_B=(\psi)^{\times p}\cdot(\mu)^{\times p}=(\psi \cdot \mu)^{\times p}=\cX(\psi\cdot\mu ; \underline{x+y})\down_B=\theta(\cR)\down_B.
\] 
Suppose $g\notin B$ and let $\rho(f;\delta)=\prod_{i=1}^{p}f_{(1)\delta^{-i+1}}\in P_{p^{k-1}}$. Using \cite[4.4.10]{JK}, we obtain 
\begin{align*}
(\varphi\cdot\lambda)(g)&=(\cX(\psi;x)\cdot\cX(\mu;y))(f;\delta)=\phi_{x}(\delta)\psi(\rho(f;\delta))\phi_{y}(\delta)\mu(\rho(f;\delta)) \\
&=\phi_{x}(\delta)\phi_{y}(\delta)\psi(\rho(f;\delta))\mu(\rho(f;\delta)) =\phi_{\underline{x+y}}(\delta)(\psi\cdot\mu)(\rho(f;\delta)) \\
&=\cX(\psi\cdot\mu;\underline{x+y})(f;\delta)=\theta(\cR)(g).
\end{align*}
This shows that $\varphi\cdot\lambda=\theta(\cR)$, so $\cT(\varphi\cdot\lambda)=\mathcal{R}$. Since $\cT(\varphi\cdot\lambda)(\emptyset)=\mathcal{R}(\emptyset)=\underline{x+y}=\underline{\cT(\varphi)(\emptyset)+\cT(\lambda)(\emptyset)}$, let $\ts\neq\emptyset$. Then from the inductive hypothesis it follows that
\begin{align*}
\cT(\varphi\cdot\lambda)(\ts)&=\mathcal{R}(\ts)=\mathcal{R}(s_1|^1\ts)=(\mathcal{R}(s_1))^{\giu}(^1\ts)=\cT(\psi\cdot\mu)(^{1}\ts) \\
&=
\begin{cases}
\cT(\psi)(^1\ts) & \ \text{if } \cT(\psi)(^1\ts)=p \\
\underline{\cT(\psi)(^1\ts) + \cT(\mu)(^1\ts)} & \ \text{if } \cT(\psi)(^1\ts)\neq p
\end{cases}
=
\begin{cases}
\cT(\varphi)(\ts) & \ \text{if } \cT(\varphi)(\ts)=p \\
\underline{\cT(\varphi)(\ts) + \cT(\lambda)(\ts)} & \ \text{if } \cT(\varphi)(\ts)\neq p
\end{cases}
\end{align*}
since $\cT(\psi)(^1\ts)=(\cT(\varphi)(s_1))^{\giu}(^1\ts)=\cT(\varphi)(s_1|^1\ts)=\cT(\varphi)(\ts)$ and, similarly, $\cT(\mu)(^1\ts)=\cT(\lambda)(\ts)$.

\noindent\underline{\textbf{\textit{Case \textnormal{$x=p$}.}}} We have that $\cT_1=\cT(\psi_1),\dotsc,\cT_p=\cT(\psi_p)$ for some mixed $\psi_1,\dotsc,\psi_p\in\Irr(P_{p^{k-1}})$, and
\[
\cT(\varphi)=(\cT(\psi_1)\mid\dots\mid\cT(\psi_p);p), \qquad \varphi=\psi_1\times\cdots\times\psi_p\up^{P_{p^k}}.
\]
Let $\mathcal{R}=(\cT(\psi_1\cdot\mu)\mid\dots\mid \cT(\psi_p\cdot\mu); p)\in\cF_{p^k}$ and $\theta(\cR)=(\psi_1\cdot\mu)\times\cdots\times (\psi_p\cdot\mu)\up^{P_{p^k}}$. From the definition of induced character (see \cite[Chapter 5, Definition 5.1]{Is}) and the fact that $B\trianglelefteq P_{p^k}$, it is straightforward to check that characters which are induced from $B$ to $P_{p^k}$ vanish when evaluated at elements in $P_{p^k}\smallsetminus B$. Therefore, if $g\in P_{p^k}\smallsetminus B$, then $\varphi\cdot\lambda(g)=\varphi(g)\lambda(g)=\varphi(g)=0=\theta(\cR)(g)$. Let $g=(f;1)=(f_1,\dotsc,f_p;1)\in B$ and consider $\gamma=(1, 2,\dotsc, p)\in P_p$. By Lemmas~\ref{lem:clif} and~\ref{lem:chP-trees}, we obtain 
\begin{align*}
(\varphi\cdot\lambda)(g) &= \varphi(f;1)\lambda(f;1) 
= (\mu)^{\times p}(f_1,\dotsc,f_p)\varphi(f_1,\dotsc,f_p;1)= \left(\prod_{i=1}^{p}\mu(f_i)\right) \left(\sum_{j=0}^{p-1} \prod_{i=1}^{p} \psi_{(i)\gamma^j}(f_i)\right) \\
&=\sum_{j=0}^{p-1} \prod_{i=1}^{p} \psi_{(i)\gamma^j}(f_i)\mu(f_i)=\sum_{j=0}^{p-1} \prod_{i=1}^{p} (\psi_{(i)\gamma^j}\cdot\mu)(f_i) \\
&=(\psi_1\cdot\mu)\times\cdots\times(\psi_p\cdot\mu)\up^{P_{p^k}}(f_1,\dotsc,f_p;1)=\theta(\cR)(g).
\end{align*}
This implies $\varphi\cdot\lambda=\theta(\cR)$, and hence $\cT(\varphi\cdot\lambda)=\mathcal{R}$. Note that $\cT(\varphi\cdot\lambda)(\emptyset)=\mathcal{R}(\emptyset)=p=\cT(\varphi)(\emptyset)$, and 
\begin{align*}
\cT(\varphi\cdot\lambda)(\ts)&=\mathcal{
R}(s_1|^1\ts)=(\mathcal{R}(s_1))^{\giu}(^{1}\ts)=\cT(\psi_{s_1}\cdot\mu)(^{1}\ts) \\
&=
\begin{cases}
\cT(\psi_{s_1})(^1\ts) & \ \text{if } \cT(\psi_{s_1})(^1\ts)=p \\
\underline{\cT(\psi_{s_1})(^1\ts) + \cT(\mu)(^1\ts)} & \ \text{if } \cT(\psi_{s_1})(^1\ts)\neq p
\end{cases}
=
\begin{cases}
\cT(\varphi)(\ts) & \ \text{if } \cT(\varphi)(\ts)=p \\
\underline{\cT(\varphi)(\ts) + \cT(\lambda)(\ts)} & \ \text{if } \cT(\varphi)(\ts)\neq p.
\end{cases}
\end{align*}
for any $\ts\in\sk_{p^k} \smallsetminus \{ \emptyset \}$, by the inductive hypothesis. This completes the proof. 
\end{proof}

\section{Equivalence relations on $\Irr(P_n)$}\label{sec:applications}\label{sec:eqonch}
This concluding part presents some consequences of Sections~\ref{sec:chT-functions} and~\ref{sec:norgalIrr(P)} on character induction from $P_n$ to $S_n$, for each prime $p$ and natural number $n$. Proposition~\ref{rem:conclusions} shows that the results on linear characters established in \cite{law,gll} and collected in Theorem~\ref{th:lawgll} do not extend to the full set of irreducible characters $\Irr(P_n)$. Before discussing the details, we fix some basic notation.

\begin{notation}
We consider the equivalence relations on $\Irr(P_n)$ induced by Galois conjugation and $N_n$-conjugation, denoted $\sim_{\cG}$ and $\sim_{N_n}$, respectively. The relation on $\Irr(P_n)$ given by having the same induction to $S_n$, written simply as $\up$, is also an equivalence relation.
\end{notation}

It is well known that if $\theta,\zeta\in\Irr(P_n)$ and either $\theta\sim_{\cG}\zeta$ or $\theta\sim_{N_n}\zeta$, then $\theta\up=\zeta\up$. Moreover, in \cite{law} it was observed that Galois conjugation implies $N_n$-conjugation among linear characters, while in \cite{gll} it was shown that the $N_n$-orbit of a linear character $\nu\in\Lin(P_n)$ is uniquely determined by the induced character $\nu\up$. We summarise these results in the following theorem.

\begin{theorem}\label{th:lawgll} \textnormal{(\cite[Subsection 4.4, p.~102]{law},\cite[Theorem A]{gll})}
Let $p$ be a prime and $n\in\N$. Given any $\lambda,\mu\in\Lin(P_n)$:
\begin{itemize}
	\item[(i)] if $\lambda\sim_{\cG}\mu$, then $\lambda\sim_{N_n}\mu$;
	\item[(ii)] $\lambda\sim_{N_n}\mu$ if and only if $\lambda\up=\mu\up$.
\end{itemize}
\end{theorem}

We now examine the limits of Theorem~\ref{th:lawgll}, showing that it fails for higher-degree characters. A counterexample is given in the proof of the following proposition. 

\begin{proposition}\label{rem:conclusions}
Theorem~\ref{th:lawgll} does not hold beyond the linear case.
\end{proposition}

\begin{proof}
We show that for some prime $p$ and $n\in\N$ there exist characters in $\Irr(P_n)$ that are Galois \mbox{conjugate} but not $N_n$-conjugate, and consequently characters whose induction from $P_n$ to $S_n$ yields the same character although they do not lie in the same $N_n$-orbit. For instance, consider $p=5$, \mbox{$n=5^3=125$} and the Sylow $5$-subgroup $P_{125}\le S_{125}$. Fix $a=2$, a primitive root modulo $5$. By Definition~\ref{def:sigmajk}, we have $\tau=(3,4,2,1)$ (see also Example~\ref{ex:111}). By Definition~\ref{def:sigma}, $\sigma\in\cG=\Gal(\Q(e^{2\pi i/5})\mid\Q)$ satisfies $(e^{2\pi i/5})^{\sigma}=(e^{2\pi i/5})^3$.
Let $t_A,t_B$ be the following $125$-labeling functions: 
\begin{itemize}
	\item[($t_A$)] $(t_A)^{-1}(0)=\{2,11,12,21,22,23,24,25,31,32,33,34,35,41,42,43,44,45,51
	,52,53,54,55\}, \\ (t_A)^{-1}(1)=\{3,13,14,15\}, \  (t_A)^{-1}(2)=\{4\}, \  (t_A)^{-1}(3)=\{5\}, \  (t_A)^{-1}(4)=\emptyset, \  (t_A)^{-1}(5)=\{ \emptyset, 1 \}$.
	\item[($t_B$)] $(t_B)^{-1}(0)=\{2,11,12,21,22,23,24,25,31,32,33,34,35,41,42,43,
	44,45,51,52,53,54,55
	\}, \\ (t_B)^{-1}(1)=\{4\}, \  (t_B)^{-1}(2)=\emptyset, \  (t_B)^{-1}(3)=\{3,13,14,15\}, \  (t_B)^{-1}(4)=\{5\},  (t_B)^{-1}(5)=\{ \emptyset, 1\}$.
\end{itemize}
For an illustration of the associated trees, see Figure~\ref{fig:counterex}. By Theorems~\ref{th:1} and~\ref{th:3}, the corresponding $\cT$-functions satisfy $\cT_B=(\cT_A)^{\sigma}$ but $\cT_B\neq (\cT_A)^{h}$ for every $h\in N_{125}$, and hence the characters $\theta=\theta(\cT_A)$ and $\theta^{\sigma}=\theta(\cT_B)$ are Galois conjugate but not $N_{125}$-conjugate. This gives the required counterexample.

\begin{figure}[htbp]
    \centering
    \begin{subfigure}[b]{0.9\textwidth}
    \centering
    \begin{tikzpicture}[scale=1.8, every node/.style={font=\small}]
        \node (R) at (0,0.2) {$\overset{5}{\bullet}$};

        \node (1) at (-3,-0.8) {$\overset{5}{\bullet}$};
        \node (2) at (-1.5,-0.8) {$\overset{0}{\bullet}$};
        \node (3) at (0,-0.8) {$\overset{1}{\bullet}$};
        \node (4) at (1.5,-0.8) {$\overset{2}{\bullet}$};
        \node (5) at (3,-0.8) {$\overset{3}{\bullet}$};        
        \draw (R) -- (1);
        \draw (R) -- (2);
        \draw (R) -- (3);
        \draw (R) -- (4);
        \draw (R) -- (5);

        \node (11) at (-3.6,-1.8) {$\overset{0}{\bullet}$};
        \node (12) at (-3.3,-1.8) {$\overset{0}{\bullet}$};
        \node (13) at (-3,-1.8) {$\overset{1}{\bullet}$};
        \node (14) at (-2.7,-1.8) {$\overset{1}{\bullet}$};
        \node (15) at (-2.4,-1.8) {$\overset{1}{\bullet}$};        
        \draw (1) -- (11);
        \draw (1) -- (12);
        \draw (1) -- (13);
        \draw (1) -- (14);
        \draw (1) -- (15);

        \node (21) at (-2.1,-1.8) {$\overset{0}{\bullet}$};
        \node (22) at (-1.8,-1.8) {$\overset{0}{\bullet}$};
        \node (23) at (-1.5,-1.8) {$\overset{0}{\bullet}$};
        \node (24) at (-1.2,-1.8) {$\overset{0}{\bullet}$};
        \node (25) at (-0.9,-1.8) {$\overset{0}{\bullet}$};
        \draw (2) -- (21);
        \draw (2) -- (22);
        \draw (2) -- (23);
        \draw (2) -- (24);
        \draw (2) -- (25);

        \node (31) at (-0.6,-1.8) {$\overset{0}{\bullet}$};
        \node (32) at (-0.3,-1.8) {$\overset{0}{\bullet}$};
        \node (33) at (0,-1.8) {$\overset{0}{\bullet}$};
        \node (34) at (0.3,-1.8) {$\overset{0}{\bullet}$};
        \node (35) at (0.6,-1.8) {$\overset{0}{\bullet}$};
        \draw (3) -- (31);
        \draw (3) -- (32);
        \draw (3) -- (33);
        \draw (3) -- (34);
        \draw (3) -- (35);

        \node (41) at (0.9,-1.8) {$\overset{0}{\bullet}$};
        \node (42) at (1.2,-1.8) {$\overset{0}{\bullet}$};
        \node (43) at (1.5,-1.8) {$\overset{0}{\bullet}$};
        \node (44) at (1.8,-1.8) {$\overset{0}{\bullet}$};
        \node (45) at (2.1,-1.8) {$\overset{0}{\bullet}$};
        \draw (4) -- (41);
        \draw (4) -- (42);
        \draw (4) -- (43);
        \draw (4) -- (44);
        \draw (4) -- (45);

        \node (51) at (2.4,-1.8) {$\overset{0}{\bullet}$};
        \node (52) at (2.7,-1.8) {$\overset{0}{\bullet}$};
        \node (53) at (3,-1.8) {$\overset{0}{\bullet}$};
        \node (54) at (3.3,-1.8) {$\overset{0}{\bullet}$};
        \node (55) at (3.6,-1.8) {$\overset{0}{\bullet}$};
        \draw (5) -- (51);
        \draw (5) -- (52);
        \draw (5) -- (53);
        \draw (5) -- (54);
        \draw (5) -- (55);
    \end{tikzpicture}
    \caption*{\scalebox{0.9}{Admissible tree associated with $t_A$, whose equivalence class corresponds to the irreducible character $\theta$.}}
    \label{fig:counterext_1}
    \end{subfigure}

\vspace{1em} 

\begin{subfigure}[b]{0.9\textwidth}
    \centering
    \begin{tikzpicture}[scale=1.8, every node/.style={font=\small}]
        \node (R) at (0,0.2) {$\overset{5}{\bullet}$};

        \node (1) at (-3,-0.8) {$\overset{5}{\bullet}$};
        \node (2) at (-1.5,-0.8) {$\overset{0}{\bullet}$};
        \node (3) at (0,-0.8) {$\overset{3}{\bullet}$};
        \node (4) at (1.5,-0.8) {$\overset{1}{\bullet}$};
        \node (5) at (3,-0.8) {$\overset{4}{\bullet}$};        
        \draw (R) -- (1);
        \draw (R) -- (2);
        \draw (R) -- (3);
        \draw (R) -- (4);
        \draw (R) -- (5);

        \node (11) at (-3.6,-1.8) {$\overset{0}{\bullet}$};
        \node (12) at (-3.3,-1.8) {$\overset{0}{\bullet}$};
        \node (13) at (-3,-1.8) {$\overset{3}{\bullet}$};
        \node (14) at (-2.7,-1.8) {$\overset{3}{\bullet}$};
        \node (15) at (-2.4,-1.8) {$\overset{3}{\bullet}$};        
        \draw (1) -- (11);
        \draw (1) -- (12);
        \draw (1) -- (13);
        \draw (1) -- (14);
        \draw (1) -- (15);

        \node (21) at (-2.1,-1.8) {$\overset{0}{\bullet}$};
        \node (22) at (-1.8,-1.8) {$\overset{0}{\bullet}$};
        \node (23) at (-1.5,-1.8) {$\overset{0}{\bullet}$};
        \node (24) at (-1.2,-1.8) {$\overset{0}{\bullet}$};
        \node (25) at (-0.9,-1.8) {$\overset{0}{\bullet}$};
        \draw (2) -- (21);
        \draw (2) -- (22);
        \draw (2) -- (23);
        \draw (2) -- (24);
        \draw (2) -- (25);

        \node (31) at (-0.6,-1.8) {$\overset{0}{\bullet}$};
        \node (32) at (-0.3,-1.8) {$\overset{0}{\bullet}$};
        \node (33) at (0,-1.8) {$\overset{0}{\bullet}$};
        \node (34) at (0.3,-1.8) {$\overset{0}{\bullet}$};
        \node (35) at (0.6,-1.8) {$\overset{0}{\bullet}$};
        \draw (3) -- (31);
        \draw (3) -- (32);
        \draw (3) -- (33);
        \draw (3) -- (34);
        \draw (3) -- (35);

        \node (41) at (0.9,-1.8) {$\overset{0}{\bullet}$};
        \node (42) at (1.2,-1.8) {$\overset{0}{\bullet}$};
        \node (43) at (1.5,-1.8) {$\overset{0}{\bullet}$};
        \node (44) at (1.8,-1.8) {$\overset{0}{\bullet}$};
        \node (45) at (2.1,-1.8) {$\overset{0}{\bullet}$};
        \draw (4) -- (41);
        \draw (4) -- (42);
        \draw (4) -- (43);
        \draw (4) -- (44);
        \draw (4) -- (45);

        \node (51) at (2.4,-1.8) {$\overset{0}{\bullet}$};
        \node (52) at (2.7,-1.8) {$\overset{0}{\bullet}$};
        \node (53) at (3,-1.8) {$\overset{0}{\bullet}$};
        \node (54) at (3.3,-1.8) {$\overset{0}{\bullet}$};
        \node (55) at (3.6,-1.8) {$\overset{0}{\bullet}$};
        \draw (5) -- (51);
        \draw (5) -- (52);
        \draw (5) -- (53);
        \draw (5) -- (54);
        \draw (5) -- (55);
    \end{tikzpicture}
    \caption*{\scalebox{0.9}{Admissible tree associated with $t_B$, whose equivalence class corresponds to the irreducible character $\theta^{\sigma}$.} }
    \label{fig:counterext_2}
    \end{subfigure}
    
    \caption{An example in $S_{125}$ involving two Galois conjugate characters in $\Irr(P_{125})$ of \mbox{degree $25$}, that are not $N_{125}$-conjugate.}
    \label{fig:counterex}
\end{figure}

\begin{figure}[htbp]
    \centering
	\begin{subfigure}[b]{0.9\textwidth}
    	\centering
    	\begin{tikzpicture}[scale=1.8, every node/.style={font=\small}]
        \node (R) at (0,0.2) {$\overset{5}{\bullet}$};

        \node (1) at (-3,-0.8) {$\overset{5}{\bullet}$};
        \node (2) at (-1.5,-0.8) {$\overset{0}{\bullet}$};
        \node (3) at (0,-0.8) {$\overset{3}{\bullet}$};
        \node (4) at (1.5,-0.8) {$\overset{1}{\bullet}$};
        \node (5) at (3,-0.8) {$\overset{4}{\bullet}$};        
        \draw (R) -- (1);
        \draw (R) -- (2);
        \draw (R) -- (3);
        \draw (R) -- (4);
        \draw (R) -- (5);

        \node (11) at (-3.6,-1.8) {$\overset{1}{\bullet}$};
        \node (12) at (-3.3,-1.8) {$\overset{0}{\bullet}$};
        \node (13) at (-3,-1.8) {$\overset{1}{\bullet}$};
        \node (14) at (-2.7,-1.8) {$\overset{0}{\bullet}$};
        \node (15) at (-2.4,-1.8) {$\overset{1}{\bullet}$};        
        \draw (1) -- (11);
        \draw (1) -- (12);
        \draw (1) -- (13);
        \draw (1) -- (14);
        \draw (1) -- (15);

        \node (21) at (-2.1,-1.8) {$\overset{0}{\bullet}$};
        \node (22) at (-1.8,-1.8) {$\overset{0}{\bullet}$};
        \node (23) at (-1.5,-1.8) {$\overset{0}{\bullet}$};
        \node (24) at (-1.2,-1.8) {$\overset{0}{\bullet}$};
        \node (25) at (-0.9,-1.8) {$\overset{0}{\bullet}$};
        \draw (2) -- (21);
        \draw (2) -- (22);
        \draw (2) -- (23);
        \draw (2) -- (24);
        \draw (2) -- (25);

       \node (31) at (-0.6,-1.8) {$\overset{0}{\bullet}$};
        \node (32) at (-0.3,-1.8) {$\overset{0}{\bullet}$};
        \node (33) at (0,-1.8) {$\overset{0}{\bullet}$};
        \node (34) at (0.3,-1.8) {$\overset{0}{\bullet}$};
        \node (35) at (0.6,-1.8) {$\overset{0}{\bullet}$};
        \draw (3) -- (31);
        \draw (3) -- (32);
        \draw (3) -- (33);
        \draw (3) -- (34);
        \draw (3) -- (35);

        \node (41) at (0.9,-1.8) {$\overset{0}{\bullet}$};
        \node (42) at (1.2,-1.8) {$\overset{0}{\bullet}$};
        \node (43) at (1.5,-1.8) {$\overset{0}{\bullet}$};
        \node (44) at (1.8,-1.8) {$\overset{0}{\bullet}$};
        \node (45) at (2.1,-1.8) {$\overset{0}{\bullet}$};
        \draw (4) -- (41);
        \draw (4) -- (42);
        \draw (4) -- (43);
        \draw (4) -- (44);
        \draw (4) -- (45);

        \node (51) at (2.4,-1.8) {$\overset{0}{\bullet}$};
        \node (52) at (2.7,-1.8) {$\overset{0}{\bullet}$};
        \node (53) at (3,-1.8) {$\overset{0}{\bullet}$};
        \node (54) at (3.3,-1.8) {$\overset{0}{\bullet}$};
        \node (55) at (3.6,-1.8) {$\overset{0}{\bullet}$};
        \draw (5) -- (51);
        \draw (5) -- (52);
        \draw (5) -- (53);
        \draw (5) -- (54);
        \draw (5) -- (55);
    \end{tikzpicture}
    \caption*{\scalebox{0.9}{Admissible tree associated with $t_C$, whose equivalence class corresponds to the irreducible character $\theta^{\sigma_2^{(3)}}$.} }
    \label{fig:op_2}
    \end{subfigure}
    
    \caption{An example in $S_{125}$ of two characters in $\Irr(P_{125})$ of degree $25$ whose induced characters from $P_{125}$ to $S_{125}$ are the same, but are not $\sim_{0,5}$-equivalent.}
    \label{fig:indno0p}
\end{figure}
\end{proof}

We care to mention that a counterexample confirming that statement (ii) of Theorem~\ref{th:lawgll} does not hold beyond the linear case was independently observed by Gabriel Navarro.  

We record a result of \cite{GL3} (see \cite[Definition 3.15 and Theorem 3.16]{GL3}), restated in Theorem~\ref{th:o,pinduction} below in terms of $\cT$-functions, which provides further information on character induction.

\begin{definition}\label{def:0,peq}\textnormal{(\cite[Definition 3.15]{GL3})}
Let $p$ be a prime, $k\geq 1$ an integer and $n$ a natural number with $p$-adic expansion as in Notation~\ref{not:p-adicn}.
\begin{itemize}
	\item[(i)] Given $\cT_1,\cT_2\in\cF_{p^k}$, we set $\cT_1\sim_{0,p}\cT_2$ if and only if there exists $t_1\in\cT_1$, $t_2\in\cT_2$ such that $(t_1)^{-1}(0)=(t_2)^{-1}(0)$ and $(t_2)^{-1}(p)=(t_2)^{-1}(p)$.
	\item[(ii)] Given $\theta,\zeta\in\Irr(P_{p^k})$, we say that $\theta\sim_{0,p}\zeta$ if and only if $\cT(\theta)\sim_{0,p}\cT(\zeta)$.
	\item[(iii)] Given any $\theta=\theta_1\times\cdots\times\theta_{q_t}, \zeta=\zeta_1\times\cdots\times\zeta_{q_t}\in\Irr(P_n)$, we say that $\theta\sim_{0,p}\zeta$ if and only if $\theta_i \sim_{0,p} \zeta_i$ for every $i\in[1,q_t]$.
\end{itemize}
\end{definition}

\begin{example}
The irreducible characters of $P_{125}$ represented in Figure~\ref{fig:counterex} are $\sim_{0,5}$-equivalent. 
\end{example}

\begin{remark}
Two irreducible characters of $P_{p^k}$ are $\sim_{0,p}$-equivalent if and only if their corresponding \mbox{$\cT$-functions} admit representatives whose preimages of $0$ and of $p$ coincide, considered separately. As shown in \cite[Definition 3.15]{GL3}, this defines an equivalence relation on $\Irr(P_{p^k})$. Intuitively, the \mbox{$\sim_{0,p}$-orbit} of an irreducible character $\theta\in\Irr(P_{p^k})$ can be described as follows. Let $t\in\cT(\theta)$ be an admissible labeling function representing the character. For each $\ts$, if $t(\ts)\in [1,p-1]$, we may freely replace $t(\ts)$ by any other value in $[1,p-1]$. If $t(\ts)\in\{0,p\}$, we leave it unchanged, so the preimage of $0$ and the preimage of $p$ are preserved. We then retain only the admissible labeling functions obtained through this process. The corresponding $\cT$-functions determine the $\sim_{0,p}$-equivalence class of $\theta$. This construction is well defined and independent of the choice of representatives, and it extends immediately to arbitrary $n$.
\end{remark}

\begin{theorem}\label{th:o,pinduction} \textnormal{(\cite[Theorem 3.16]{GL3})}
Let $p$ be a prime and $n\in\N$. The $\sim_{0,p}$-equivalence is an equivalence relation on $\Irr(P_n)$. Moreover, given any $\theta,\zeta\in\Irr(P_n)$, if $\theta\sim_{0,p}\zeta$ then $\theta\up=\zeta\up$.
\end{theorem}

A natural question is whether the converse of Theorem~\ref{th:o,pinduction} holds. The following proposition shows that this is not the case. 

\begin{proposition}\label{prop:no}
The $\sim_{0,p}$-equivalence does not coincide with the property of yielding same induced character from $P_n$ to $S_n$. 
\end{proposition}

\begin{proof}
We show that for some prime $p$ and $n\in\N$ there exist characters in $\Irr(P_n)$ whose induction from $P_n$ to $S_n$ is the same even though they do not lie in the same $\sim_{0,p}$-equivalence class. For instance, consider $p=5$, $n=5^3=125$ and the corresponding Sylow $5$-subgroup $P_{125}\le S_{125}$. Fix $a=2$, a primitive root modulo $5$, and define $\tau=(3,4,2,1)$ (see Example~\ref{ex:111}). Let $t_A$ be the $125$-labeling function as in the proof of Proposition~\ref{rem:conclusions}, and let $t_C$ be the $125$-labeling function defined by 
\begin{itemize}
	\item[] $(t_C)^{-1}(0)=\{2,12,14,21,22,23,24,25,31,32,33,34,35,
41,42,43,44,45,51,52,53,54,55\}, \\ (t_C)^{-1}(1)=\{4,11,13,15\}, \  (t_C)^{-1}(2)=\emptyset, \  (t_C)^{-1}(3)=\{3\}, \  (t_C)^{-1}(4)=\{5\}, \  (t_C)^{-1}(5)=\{ \emptyset, 1 \}$.
\end{itemize} 
We refer the reader to Figures~\ref{fig:counterex} and ~\ref{fig:indno0p} to visualize the associated trees. The corresponding $\cT$-functions \mbox{satisfy} $\cT_C=(\cT_A)^{\sigma_2^{(3)}}$. We now prove that $\cT_A\not\sim_{0,5}\cT_C$, arguing by contradiction. Suppose that \mbox{$\cT_A\sim_{0,5}\cT_C$}. Using Notation~\ref{not:T=()}, there exist some $r=(r_1\mid\dots\mid r_5;5)\in\cT_A$ and \mbox{$l=(l_1\mid\dots\mid l_5;5)\in\cT_C$} such that $r^{-1}(0)=l^{-1}(0)$ and $r^{-1}(5)=l^{-1}(5)$. In particular, this implies $r_i^{-1}(0)=l_i^{-1}(0)$ and \mbox{$r_i^{-1}(5)=l_i^{-1}(5)$} for all $i\in [1,5]$. Therefore, with $t_A$ and $t_C$ as above, there exists some $j\in [1,5]$ such that \mbox{\scalebox{0.95}{$(\cT_A(1))^{\giu}\sim_{0,5}(\cT_C(j))^{\giu}$}}. If $j\neq 1$, then $((t_C(j))^{\giu})^{-1}(5)=\emptyset\neq ((t_A(1))^{\giu})^{-1}(5)$, and the corresponding $\cT$-functions cannot be \mbox{$\sim_{0,5}$-equivalent}. Hence $j=1$ and $(\cT_A(1))^{\giu}\sim_{0,5}(\cT_C(1))^{\giu}$. By Definition~\ref{def:0,peq}, there exist $f\in (\cT_A(1))^{\giu}$ and $g\in(\cT_C(1))^{\giu}$ such that $f^{-1}(0)=g^{-1}(0)$. Since $f$ represents $(\cT_A(1))^{\giu}$, it follows that $f^{-1}(0)=\{u,v\}$ for some $u,v\in [1,5]$ such that $u-v\equiv\pm 1 \pmod{5}$. Similarly, from $g\in(\cT_C(1))^{\giu}$, we deduce that $g^{-1}(0)=\{w,z\}$ for some $w,z\in [1,5]$ such that $w-z\equiv\pm 2 \pmod{5}$. Then the \mbox{condition} $f^{-1}(0)=g^{-1}(0)$ is a contradiction, since we cannot find two numbers in $[1,5]$ whose difference is simultaneously $\pm 1$ and $\pm 2$ modulo $5$. We conclude that the irreducible characters associated to $\cT_A$ and $\cT_C$ are $N_{125}$-conjugate (and hence have same induction from $P_{125}$ to $S_{125}$) yet fail to be $\sim_{0,5}$-equivalent. This provides the desired counterexample.
\end{proof}

Finally, we compare the equivalence relations on $\Irr(P_n)$ introduced in this section and summarize their relationships in the following table.

\begin{table}[htbp]
\centering
\renewcommand{\arraystretch}{1.3}
\begin{tabular}{|l|c|c|}
\hline
\multicolumn{3}{|c|}{Let $p$ be a prime, $n\in\N$, and $\theta,\zeta\in\Irr(P_n)$.} \\
\multicolumn{3}{|c|}{If $\theta\sim\zeta$, where $\sim$ is one of $\sim_{\cG},\sim_{N_n},\sim_{0,p}$, then $\theta\up=\zeta\up$, but the converse doesn't hold.} \\ 
\multicolumn{3}{|c|}{In general, we have:} \\ \hline
\underline{Statement} & \underline{Linear case} & \underline{Non-linear case} \\ \hline
$\theta\sim_{\cG}\zeta\Rightarrow\theta\sim_{0,p}\zeta$ & true (the converse is false) & true (the converse is false) \\ \hline
$\theta\sim_{\cG}\zeta\Rightarrow\theta\sim_{N_n}\zeta$ & true (the converse is false)  & false (the converse is false) \\ \hline
$\theta\sim_{N_n}\zeta \Leftrightarrow \theta\up=\zeta\up$ & true & false ($\theta\sim_{N_n}\zeta\not\Leftarrow\theta\up=\zeta\up$) \\ \hline
$\theta\sim_{0,p}\zeta \Rightarrow \theta\sim_{N_n}\zeta$ & true (the converse is false) & false (the converse is false) \\ \hline
\multicolumn{3}{|c|}{If $n=p^k$ for some $k\in\N$,} \\ \hline
$\theta\sim_{0,p}\zeta \Leftrightarrow \theta\sim_{N_n}\zeta$ & true & false (both directions fail) \\ \hline
\end{tabular}
\caption{Comparison of the equivalence relations $\sim_{\cG},\sim_{N_n},\sim_{0,p},\up=\up$ on $\Irr(P_n)$ in the linear and non-linear cases.}
\label{tab:schema}
\end{table}

\end{document}